\documentclass[11pt,a4paper]{article}
\usepackage{mathrsfs}
\usepackage{amssymb}
\usepackage{amsfonts}
\usepackage{color,xcolor}
\usepackage{graphicx}
\usepackage{manfnt}
\RequirePackage[colorlinks,citecolor=blue,urlcolor=blue]{hyperref}
\newtheorem{theorem}{Theorem}[section]
\newtheorem{corollary}{Corollary}[section]
\newtheorem{lemma}{Lemma}[section]
\newtheorem{proposition}{Proposition}[section]

\newtheorem{example}{Example}[section]

\def\[{{\Big[}}\def\]{{\Big]}}\def\({{\Big(}}\def\){{\Big)}}

\def\cM{{\mathcal M}}

\def\mN{{\mathbb N}}

\def\={&\!\!=\!\!&}
\def\geq{\geqslant}\def\leq{\leqslant}

\begin{document}
\title{\bf On weak solutions of stochastic differential equations with sharp drift coefficients}

\author{Jinlong Wei$^a$, Guangying Lv$^b$ and Jiang-Lun Wu$^{c}$\thanks{Corresponding author. Email: j.l.wu@swansea.ac.uk}\\
\\ {\small \it $^a$ School of Statistics and Mathematics, Zhongnan University of}\\  {\small \it Economics and Law, Wuhan 430073, China}
\\ {\small \it $^b$ Institute of Applied Mathematics, Henan University} \\ {\small \it Kaifeng 475001, China}
\\ {\small \it $^c$ Department of Mathematics, Swansea University, Swansea SA2 8PP, UK}}

\date{}
\maketitle

\noindent{\hrulefill}
\vskip1mm\noindent{\bf Abstract} We extend Krylov and R\"{o}ckner's result \cite{KR} to the drift coefficients in critical Lebesgue space, and prove the existence and
uniqueness of weak solutions for a class of SDEs. To be more precise, let $b: [0,T]\times{\mathbb R}^d\rightarrow{\mathbb R}^d$ be Borel measurable, where $T>0$
is arbitrarily fixed. Consider
$$X_t=x+\int_0^tb(s,X_s)ds+W_t,\quad t\in[0,T], \, x\in{\mathbb R}^d,$$
where $\{W_t\}_{t\in[0,T]}$ is a $d$-dimensional standard Wiener process. If $b=b_1+b_2$ such that
$b_1(T-\cdot)\in\mathcal{C}_q^0((0,T];L^p({\mathbb R}^d))$ with $2/q+d/p=1$ for $p,q\ge1$ and $\|b_1(T-\cdot)\|_{\mathcal{C}_q((0,T];L^p({\mathbb R}^d))}$ is
sufficiently small, and that $b_2$ is bounded and Borel measurable, then there exits a unique weak solution to the above equation.
Furthermore, we obtain the strong Feller property of the semi-group and existence of density associated with above SDE.
Besides, we extend the classical partial differential equations (PDEs) results for $L^q(0,T;L^p({\mathbb R}^d))$ coefficients to $L^\infty_q(0,T;L^p({\mathbb R}^d))$ ones,
and derive the Lipschitz regularity for solutions of second order parabolic PDEs (see Lemma \ref{lem2.1}).

\vskip1mm\noindent {\bf MSC (2010):} 60H10; 34F05

\vskip1mm\noindent {\bf Keywords:} Existence, uniqueness, weak solution, SDEs with irregular drifts, strong Feller property

\vskip1mm\noindent{\hrulefill}

\section{Introduction}\label{sec1}
\setcounter{equation}{0}
Let $T\in (0,\infty)$ be arbitrarily fixed. For a Borel measure function $h:[0,T]\to\mathbb{R}$, we set the notation $\mathcal{I}_Th(t):=h(T-t), t\in[0,T]$. Furthermore,
for a (joint) Borel measurable function  $f: [0,T]\times{\mathbb R}^d\to{\mathbb R}$, we denote $\mathcal{I}_Tf(t,x):=f(T-t,x), (t,x)\in[0,T]\times{\mathbb R}^d$. We are concerned with the following
stochastic differential equation (SDE for short) in ${\mathbb R}^d$:
\begin{eqnarray}\label{1.1}
\left\{\begin{array}{ll}
dX_t(x)=b(t,X_t(x))dt+dW_t, \ 0<t\leq T,\\
X_0(x)=x\in{\mathbb R}^d,  \end{array}\right.
\end{eqnarray}
where $\{W_t\}_{0\leq t\leq T}=\{(W_{1,t}, W_{2,t}, _{\cdots}, W_{d,t})\}_{0\leq t\leq T}$ is a $d$-dimensional standard Wiener process defined on a given
stochastic basis ($\Omega, {\mathcal F},\{{\mathcal F}_{t}\}_{0\leq t\leq T},{\mathbb P}$), and the drift coefficient $b: [0,T]\times{\mathbb R}^d\rightarrow{\mathbb R}^d$
is Borel measurable such that $b\in L^1(0,T;L^1_{loc}({\mathbb R}^d;{\mathbb R}^d))$.

When $b$ is Lipschitz continuous in $x\in\mathbb{R}^d$ uniformly for $t\in[0,T]$, the existence and uniqueness for strong solutions of (\ref{1.1}) can be followed by the
classical Cauchy-Lipschitz theorem. This result was firstly extended by Veretennikov \cite{Ver} to bounded Borel measurable function $b$. Since then,
Veretennikov's result was strengthened in different forms under the same assumption on $b$. For instance, Mohammed, Nilssen, Proske in \cite{MNP}
not only showed the existence and uniqueness of strong solutions, but also obtained that the unique strong solution forms a Sobolev differentiable stochastic flows;
Davie showed in \cite{Dav07}  that for almost every Wiener path $W$, there is a unique continuous $X$ satisfying the integral equation (also see \cite{Fla11}).

For integrable drift coefficient, i.e.
\begin{eqnarray}\label{1.2}
b\in L^q(0,T;L^p({\mathbb R}^d;{\mathbb R}^d))
\end{eqnarray}
with some $p,q\in [2,\infty)$ such that
\begin{eqnarray}\label{1.3}
\frac{2}{q}+\frac{d}{p}<1,
\end{eqnarray}
by applying Girsanov's theorem and Krylov's estimate, Krylov and R\"{o}ckner \cite{KR} showed the existence and uniqueness of strong solutions for SDE
(\ref{1.1}). On the other hand, under the same conditions (\ref{1.2}) and (\ref{1.3}),  Fedrizzi and Flandoli \cite{FF13} proved the $\alpha$-H\"{o}lder
continuity of $x\in\mathbb{R}^d\mapsto X_t(x)\in\mathbb{R}^d$ for every $\alpha\in (0,1)$. Some further interesting extensions for non-constant diffusion
coefficients, the reader is referred to Zhang \cite{Zha05,Zha11}.

However, to the best of our knowledge, there are few investigations to consider the critical case, i.e.
\begin{eqnarray}\label{1.4}
\frac{2}{q}+\frac{d}{p}=1.
\end{eqnarray}
Inspired by Ambrosio \cite{Amb04}, by introducing a notation of Lagrangian flow, Beck, Flandoli, Gubinelliz and Maurellix in \cite{BFGM} derived the existence and
uniqueness, in the present setting, for SDE (\ref{1.1}) for every $\omega\in \Omega$ being fixed. As stated in \cite{BFGM}: ``we do not know whether or not strong
solutions exist and are unique under the conditions (\ref{1.2}) and (\ref{1.4}), and as we know this problem is still open". In fact, under the conditions (\ref{1.2}) and
(\ref{1.4}), there are as well few works for the existence and uniqueness of weak solutions to (\ref{1.1}) and yet it remains to be a challenging problem.

The above problems are the main driving source for us to work out the present paper. In this paper, we will discuss the existence and uniqueness of weak solutions
to SDE (\ref{1.1}) under the critical case (\ref{1.4}) by replacing the integrable condition (\ref{1.2}) to that $b\in \mathcal{C}_q((0,T];L^p({\mathbb R}^d;{\mathbb R}^d))$.

Our main result is the following

\begin{theorem} \label{the1.1}  Suppose that $p,q\in [1,\infty)$. Let $b=b_1+b_2$ such that $\mathcal{I}_Tb_1 \in \mathcal{C}_q((0,T];L^p({\mathbb R}^d;{\mathbb R}^d))$ with $p,q$ satisfying (\ref{1.4}),
$b_2$ is bounded. Suppose that $\|\mathcal{I}_Tb_1\|_{\mathcal{C}_q((0,T];L^p({\mathbb R}^d))}$ is sufficiently small, then we have the following consequences

(i) There is a filtered probability space $(\tilde{\Omega},\tilde{\mathcal{F}},\{\tilde{\mathcal{F}}_t\}_{0\leq t\leq T},\tilde{{\mathbb P}})$ on which there are two  processes $\{\tilde{X}_t\}_{t\in[0,T]}$
and $\{\tilde{W}_t\}_{t\in[0,T]}$ such that \{$\tilde{W}_t\}$ is a $d$-dimensional $\{\tilde{\mathcal{F}}_t\}$-Wiener process and $\{\tilde{X}_t\}$ is an $\{\tilde{\mathcal{F}}_t\}$-adapted,
continuous, $d$-dimensional process for which
\begin{eqnarray}\label{1.5}
\tilde{{\mathbb P}}\Big(\int_0^T|b(t,\tilde{X}_t)|dt<\infty\Big)=1
\end{eqnarray}
and the following equation holds
\begin{eqnarray}\label{1.6}
\tilde{X}_t=x+\int_0^tb(s,\tilde{X}_s)ds+\tilde{W}_t, \quad t\in[0,T],\,\, \tilde{{\mathbb P}}-a.s..
\end{eqnarray}

(ii) If in addition $\mathcal{I}_Tb_1\in\mathcal{C}^0_q((0,T];L^p({\mathbb R}^d;{\mathbb R}^d))$, then all weak solutions for SDE (\ref{1.1}) have the same probability law on $d$-dimensional classical
Wiener space $(W^d([0,T]),\mathcal{B}(W^d([0,T])))$. We then use ${\mathbb P}_x$ to denote the unique probability law on $(W^d([0,T]),\mathcal{B}(W^d([0,T])))$ corresponding to the initial
value $x\in\mathbb{R}^d$.

(iii) With the same condition of (ii). For every $f\in L^\infty({\mathbb R}^d)$, we define
$$P_tf(x):={\mathbb E}^{{\mathbb P}_x}f(w(t)), t>0,$$
where $w(t)$ is the canonical realisation of a weak solution $\{X_t\}$ with initial $X_0=x\in\mathbb{R}^d$
on $(W^d([0,T]),\mathcal{B}(W^d([0,T])))$. Then, the semigroup $\{P_t\}$ has strong Feller property, i.e. each $P_t$ maps a
bounded function to a bounded and continuous function. Moreover, $P_t$ admits a density $p(t,x,y)$ for almost all $t\in  [0,T]$. Besides, for every $t_0>0$ and for every $r\in [1,\infty)$,
\begin{eqnarray}\label{1.7}
\int_{t_0}^T\int_{{\mathbb R}^d}|p(t,x,y)|^rdydt<\infty.
\end{eqnarray}.
\end{theorem}

\noindent
\textbf{Remark 1.1.}  Here the existence is only for weak solutions, and the uniqueness holds only in the sense of probability laws.
However, we do not know in the present setting for general $d$, whether the strong solutions do exist and further, if they do, whether the uniqueness holds.
But for $d=1$, with the aid of Tanack, Tsuchiya, Watanabe's result \cite[Proposition 1.1]{YTW}, we can give a positive answer on strong uniqueness. To be more precise, we have
\begin{corollary} \label{cor1.1} Suppose that $p,q\in [1,\infty)$. Let $b=b_1+b_2$ such that $\mathcal{I}_Tb_1 \in \mathcal{C}_q^0((0,T];L^p({\mathbb R}))$ with $2/q+1/p=1$,
$b_2$ is bounded. Suppose that $\|\mathcal{I}_Tb_1\|_{\mathcal{C}_q((0,T];L^p({\mathbb R}^d))}$ is sufficiently small, then SDE (\ref{1.1}) exists a unique strong solution with $d=1$.
\end{corollary}

To illustrate the present result and to compare it with the results given by the case of $b\in L^q(0,T;L^p({\mathbb R}^d;{\mathbb R}^d))$, let us look at an example.

\begin{example}\label{exa1.1} Let $p,q\in [1,\infty)$ such that (\ref{1.4}) holds. Suppose $\tilde{b}\in\mathcal{C}([0,T];L^p({\mathbb R}^d;{\mathbb R}^d))$ with $T=1/2$, so that $\|\tilde{b}\|_{\mathcal{C}([0,T];L^p({\mathbb R}^d))}$ is small enough. We set
\begin{eqnarray}\label{1.8}
b(t,x)=(\frac12-t)^{-\frac{1}{q}}|\log (\frac12-t)|^{-1}\tilde{b}(t,x),
\end{eqnarray}
then
\begin{eqnarray*}
t^{\frac1q}b(\frac12-t,x)=|\log t|^{-1}\tilde{b}(\frac12-t,x),
\end{eqnarray*}
which indicates that $\mathcal{I}_{\frac12}b\in\mathcal{C}^0_q((0,\frac12];L^p({\mathbb R}^d;{\mathbb R}^d))$. By Theorem \ref{the1.1}, SDE (\ref{1.1}) with $b$ given by (\ref{1.8}) exists a unique weak solution.
On the other hand, from the explicit form (\ref{1.8}), $b\in L^q(0,\frac12;L^p({\mathbb R}^d;{\mathbb R}^d))$, but for every $q_1>q$, $b\notin L^{q_1}(0,\frac12;L^p({\mathbb R}^d;{\mathbb R}^d))$. Now since
$2/q+d/p=1$, from existing results, we do not know whether there exists a unique solution to (\ref{1.1}). From this point of view, it is clear that we extend the existing results on
$L^q(0,\frac12;L^p({\mathbb R}^d;{\mathbb R}^d))$ with $2/q+d/p<1$ to $2/q+d/p=1$ partially.
\end{example}

The rest of this paper is arranged as follows. In Section 2, we present some preliminaries. Section 3 is devoted to the proof of existence for Theorem \ref{the1.1} and in Section 4,
the final section, we prove the uniqueness, strong Feller property as well as the existence of the density.

When there is no ambiguity, we use $C$ to denote a constant whose true value may vary from line to line.  As usual, $\mN$ stands for the set of all natural numbers.

\section{Preliminaries}\label{sec2}
\setcounter{equation}{0}
To start with, let us introduce some spaces. For $q\geq 1$, we denote $L^\infty_q(0,T)$ the space of all Borel-measurable functions $h: [0,T]\rightarrow {\mathbb R}$ such that
$\sup\limits_{0\leq t\leq T}(|h(t)|t^{\frac{1}{q}})<\infty$, the norm is defined as follows
\begin{eqnarray*}
\|h\|_{L^\infty_q(0,T)}:=\sup_{0\leq t\leq T}(|h(t)|t^{\frac{1}{q}}).
\end{eqnarray*}
Clearly, $L^\infty_q(0,T)$ is a Banach space. If $h\in L^\infty_q(0,T)$, then there is a positive constant $C$ such that $|h(t)|\leq C t^{-\frac{1}{q}}$. Observing that
$t^{-\frac{1}{q}}\in L^{q,\infty}(0,T)$ (weak-$L^q(0,T)$ space), then $h\in L^{q,\infty}(0,T)$. Similarly, let $p\geq 1$, we define $L^\infty_q(0,T;L^p({\mathbb R}^d))$
to be the set of all $L^p({\mathbb R}^d)$-valued $L^\infty_q(0,T)$ functions $f$, such that
\begin{eqnarray}\label{2.1}
\|f\|_{L^\infty_q(0,T;L^p({\mathbb R}^d))}:=\sup_{0\leq t\leq T}(\|f(t)\|_{L^p({\mathbb R}^d)}t^{\frac{1}{q}})<\infty.
\end{eqnarray}
Then $L^\infty_q(0,T;L^p({\mathbb R}^d))$ is a Banach space as well.

Analogously, we define $\mathcal{C}_q((0,T])$ the space of all continuous function $h: (0,T]\rightarrow {\mathbb R}$ such that $\sup\limits_{0< t\leq T}(|h(t)|t^{\frac{1}{q}})<\infty$, the norm
is the same as on $L^\infty_q(0,T)$. $\mathcal{C}^0_q((0,T])$ is the space consisting by all the functions $h$ in $\mathcal{C}_q((0,T])$ such that $\lim_{t\downarrow 0}(|h(t)|t^{\frac{1}{q}})=0$. Then,
with the norm on $\mathcal{C}_q((0,T])$, $\mathcal{C}_q^0((0,T])$ is still a Banach space. Respectively, for $p\geq 1$, we define $\mathcal{C}_q((0,T];L^p({\mathbb R}^d))$ and $\mathcal{C}^0_q((0,T];L^p({\mathbb R}^d))$,
and the norms are given by (\ref{2.1}). Now let us give a good approximation property for functions in $\mathcal{C}_q((0,T];L^p({\mathbb R}^d))$.

\begin{proposition} \label{pro2.1} Suppose that $p,q\in [1,\infty)$. Given a function $f$ in $\mathcal{C}_q((0,T];L^p({\mathbb R}^d))$, we set $f_n(t,x)=(f(t,\cdot)\ast\rho_n)(x), n\in\mathbb{N}$, where
$\ast$ stands for the usual convolution and
\begin{eqnarray}\label{2.2}
\rho_n :=n^d \rho(n\cdot)  \ with \  0\leq \rho \in \mathcal{C}^\infty_0({\mathbb R}^d) , \ \ support(\rho)\subset B_0(1), \ \int_{{\mathbb R}^d}\rho(x)dx=1.
\end{eqnarray}
Then
\begin{eqnarray}\label{2.3}
\lim_{n\rightarrow \infty}\sup_{0<t\leq T}\Big(t^{\frac1q}\|f_n(t)-f(t)\|_{L^p({\mathbb R}^d)}\Big)=0.
\end{eqnarray}
\end{proposition}
\vskip2mm\noindent
\textbf{Proof.} Notice that $f\in \mathcal{C}_q((0,T];L^p({\mathbb R}^d))$  means that $t^{\frac1q}f\in \mathcal{C}([0,T];L^p({\mathbb R}^d))$, if one defines the value at $0$ by its right limit. Thus,
to prove (\ref{2.3}), it is sufficient to show that for $f\in \mathcal{C}([0,T];L^p({\mathbb R}^d))$
\begin{eqnarray}\label{2.4}
\limsup_{n\rightarrow \infty}\sup_{0\leq t\leq T}\|f_n(t)-f(t)\|_{L^p({\mathbb R}^d)}=0.
\end{eqnarray}
By virtue of properties of the convolution, for every fixed $t\in [0,T]$, then
\begin{eqnarray}\label{2.5}
\lim_{n\rightarrow \infty}\|f_n(t)-f(t)\|_{L^p({\mathbb R}^d)}=0.
\end{eqnarray}

On the other hand for $t_1,t_2\in[0,T]$, by utilising Young's inequality,
\begin{eqnarray}\label{2.6}
\|f_n(t_1)-f_n(t_2)\|_{L^p({\mathbb R}^d)}^p&=&\int_{{\mathbb R}^d}|(f(t_1,\cdot)-f(t_2,\cdot))
\ast\rho_n(x)|^pdx\cr\cr&\leq& \int_{{\mathbb R}^d}|f(t_1,x)-f(t_2,x)|^pdx.
\end{eqnarray}
From (\ref{2.6}), for any $\epsilon>0$, there exists $\delta>0$ such that for $|t_1-t_2|\leq\delta$, then one has uniformly in $n$ the following
\begin{eqnarray}\label{2.7}
\|f_n(t_1)-f_n(t_2)\|_{L^p({\mathbb R}^d)}\leq \|f(t_1)-f(t_2)\|_{L^p({\mathbb R}^d)}<\frac\epsilon2.
\end{eqnarray}

Let $t\in [0,T]$ be given, then (\ref{2.5}) holds. With the aid of (\ref{2.6}) and (\ref{2.7}), then
\begin{eqnarray*}
&&\limsup_{n\rightarrow \infty}\sup_{[t-\delta,t+\delta]\cap [0,T]}\|f_n(s)-f(s)\|_{L^p({\mathbb R}^d)}\cr\cr&\leq& \limsup_{n\rightarrow \infty}
\sup_{[t-\delta,t+\delta]\cap [0,T]}\|f_n(s)-f(s)-f_n(t)+f(t)\|_{L^p({\mathbb R}^d)} +\limsup_{n\rightarrow \infty}\|f_n(t)-f(t)\|_{L^p({\mathbb R}^d)}\cr\cr&<&\epsilon.
\end{eqnarray*}
Since $\epsilon>0$ and $t\in [0,T]$ are arbitrary, we conclude that (\ref{2.4}) holds. $\Box$

\vskip4mm\noindent
\textbf{Remark 2.1.} We claim that the above approximating property is not true if one takes the function in $L^\infty_q((0,T];L^p({\mathbb R}^d))$ instead of in
$\mathcal{C}_q((0,T];L^p({\mathbb R}^d))$. To show this, one needs to claims that (\ref{2.4}) is not true if $f\in L^\infty(0,T;L^p({\mathbb R}^d))$. We give a counter example
below. For simplicity, we assume that $T=d=1$ and $p=2$. For $k\in\mN$, we define $f_k(x)$ by the following
\begin{eqnarray*}
f^k(x):=k1_{[k,k+\frac{1}{k^2})}(x),
\end{eqnarray*}
and further set
\begin{eqnarray*}
f(t,x):=\sum_{k=1}^\infty 1_{[\frac{k-1}{k},\frac{k}{k+1})}(t)f^k(x)=\sum_{k=1}^\infty 1_{[\frac{k-1}{k},\frac{k}{k+1})}(t)k1_{[k,k+\frac{1}{k^2})}(x).
\end{eqnarray*}
Then
\begin{eqnarray*}
\int_{{\mathbb R}}|f(t,x)|^2dx=\sum_{k=1}^\infty 1_{[\frac{k-1}{k},\frac{k}{k+1})}(t)\int_{{\mathbb R}}k^21_{[k,k+\frac{1}{k^2})}(x)dx
=\sum_{k=1}^\infty 1_{[\frac{k-1}{k},\frac{k}{k+1})}(t)=1_{[0,1)}(t).
\end{eqnarray*}
Hence $f\in L^\infty(0,1;L^2({\mathbb R}))$. We estimate (\ref{2.4}) by the following
\begin{eqnarray*}
&&\int_{{\mathbb R}}|f_n(t,x)-f(t,x)|^2dx\cr\cr&=& \int_{{\mathbb R}}\Big|\int_{{\mathbb R}}f(t,y)\rho_n(x-y)dy-f(t,x)\Big|^2dx
\cr\cr&=&\int_{{\mathbb R}}\sum_{k=1}^\infty k^2 1_{[\frac{k-1}{k},\frac{k}{k+1})}(t)\Big|\int_{{\mathbb R}}1_{[k,k+\frac{1}{k^2})}(y)\rho_n(x-y)dy
-1_{[k,k+\frac{1}{k^2})}(x)\Big|^2dx\cr\cr&=&\int_{{\mathbb R}}\sum_{k=1}^\infty k^2 1_{[\frac{k-1}{k},\frac{k}{k+1})}(t)\Big|\int_k^{k+\frac{1}{k^2}}\rho_n(x-y)dy
-1_{[k,k+\frac{1}{k^2})}(x)\Big|^2dx\cr\cr&\geq& 1_{[\frac{k-1}{k},\frac{k}{k+1})}(t)\int_k^{k+\frac{1}{k^2}}k^2 \Big|\int_k^{k+\frac{1}{k^2}}\rho_n(x-y)dy
-1\Big|^2dx.
\end{eqnarray*}
For any fixed $n$, for sufficiently large $k$, we have $|\int_k^{k+\frac{1}{k^2}}\rho_n(x-y)dy|<\frac{1}{2}$. Thus
\begin{eqnarray*}
&&\sup_{0\leq t\leq 1}\int_{{\mathbb R}}|f_n(t,x)-f(t,x)|^2dx
\cr\cr&\geq& \sup_{k}\sup_{0\leq t\leq 1}1_{[\frac{k-1}{k},\frac{k}{k+1})}(t)\int_k^{k+\frac{1}{k^2}}k^2 \Big|\int_k^{k+\frac{1}{k^2}}\rho_n(x-y)dy
-1\Big|^2dx
\geq\frac{1}{4}.
\end{eqnarray*}
Therefore
\begin{eqnarray*}
\liminf_{n\rightarrow\infty}\sup_{0\leq t\leq 1}\int_{{\mathbb R}}|f_n(t,x)-f(t,x)|^2dx\geq\frac{1}{4}.
\end{eqnarray*}

\vskip2mm\par
Let $p\geq1$ such that $g\in L^1(0,T;L^p_{loc}({\mathbb R}^d;{\mathbb R}^d))$, $f\in L^1(0,T;L^1_{loc}({\mathbb R}^d))$. Consider the following Cauchy problem for $u:[0,T]\times\mathbb{R}^d\to\mathbb{R}$
\begin{eqnarray}\label{2.8}
\left\{\begin{array}{ll}
\partial_{t}u(t,x)=\frac{1}{2}\Delta u(t,x)+g(t,x)\cdot \nabla u(t,x)+f(t,x), \ (t,x)\in (0,T)\times {\mathbb R}^d, \\
u(0,x)=0, \  x\in{\mathbb R}^d.  \end{array}\right.
\end{eqnarray}
We call $u(t,x)$ a generalized solution of (\ref{2.8}) if it lies in $\mathcal{C}([0,T];W^{1,p}({\mathbb R}^d))$ such that for every test function $\varphi\in \mathcal{C}_0^\infty([0,T)\times {\mathbb R}^d)$, the following holds
\begin{eqnarray}\label{2.9}
0&=&\int_0^T\int_{{\mathbb R}^d}u(t,x)\partial_t\varphi(t,x)dxdt+\frac{1}{2}\int_0^T
\int_{{\mathbb R}^d}u(t,x)\Delta \varphi(t,x)dxdt\cr\cr&&+\int_0^T
\int_{{\mathbb R}^d}g(t,x)\cdot \nabla u(t,x) \varphi(t,x)dxdt+\int_0^T
\int_{{\mathbb R}^d}f(t,x)\varphi(t,x)dxdt.
\end{eqnarray}

The following proposition is routine and we therefore omit its proof. For more details, the reader is referred to \cite[Proposition 3.5]{Zha13}.
\begin{proposition} \label{pro2.2} Let $p\in[1,\infty)$ such that $g\in L^1(0,T;L^p({\mathbb R}^d;{\mathbb R}^d))$, $f\in L^1(0,T;L^p({\mathbb R}^d))$ and $u\in\mathcal{C}([0,T];W^{1,p}({\mathbb R}^d))$. The following statements
are equivalent

(i) $u$ is a generalized solution of (\ref{2.8}).

(ii) For every $\psi\in \mathcal{C}_0^\infty({\mathbb R}^d)$, and every $t\in [0,T)$, the following holds
\begin{eqnarray*}
\int_{{\mathbb R}^d}u(t,x)\psi(x)dx&=&\frac{1}{2}\int_0^t
\int_{{\mathbb R}^d}u(s,x)\Delta \psi(x)dxds+\int_0^t\int_{{\mathbb R}^d}g(s,x)\cdot \nabla u(s,x) \psi(x)dxds\cr\cr&&+\int_0^t
\int_{{\mathbb R}^d}f(s,x)\psi(x)dxds.
\end{eqnarray*}

(iii) For every $t\in [0,T]$ and for almost everywhere $x\in{\mathbb R}^d$, $u$ fulfils the following integral equation
\begin{eqnarray}\label{2.10}
u(t,x)=\int_0^tK(t-s,\cdot)\ast (g(s,\cdot)\cdot \nabla u(s,\cdot))(x)ds+\int_0^t(K(t-s,\cdot)\ast f(s,\cdot))(x)ds
\end{eqnarray}
where $K(t,x)=(2\pi t)^{-\frac{d}{2}}e^{-\frac{|x|^2}{2t}}, t>0, x\in\mathbb{R}^d$.
\end{proposition}

We now state a useful lemma.
\begin{lemma} \label{lem2.1}  Let $p,q\in [1,\infty)$ and $g\in L^\infty_q(0,T;L^p({\mathbb R}^d;{\mathbb R}^d))$, $f\in L^\infty_q(0,T;L^p({\mathbb R}^d))$, such that (\ref{1.4}) holds true and
$\|g\|_{L^\infty_q(0,T;L^p({\mathbb R}^d))}$ is sufficiently small. Then the Cauchy problem (\ref{2.8}) has a unique generalized solution $u$. Moreover, the unique generalized
solution lies in $L^\infty(0,T;W^{1,\infty}({\mathbb R}^d))$ and there is a constant $C_0(p,d)$ such that
\begin{eqnarray}\label{2.11}
 \|u\|_{L^\infty(0,T;W^{1,\infty}({\mathbb R}^d))} \leq \frac{C_0(p,d)\|f\|_{L^\infty_q(0,T;L^p({\mathbb R}^d))}}{1-C_0(p,d)
\|g\|_{L^\infty_q(0,T;L^p({\mathbb R}^d))}}.
\end{eqnarray}
\end{lemma}

\vskip2mm\noindent\textbf{Proof.} We prove the result by first assuming that $g=0$. With the help of Proposition \ref{pro2.2}, it suffices to show that
\begin{eqnarray}\label{2.12}
u(t,x)=\int_0^{t}(K(t-s,\cdot)\ast f(s,\cdot))(x)ds
\end{eqnarray}
is in $\mathcal{C}([0,T];W^{1,p}({\mathbb R}^d))\cap L^\infty(0,T;W^{1,\infty}({\mathbb R}^d))$. Firstly, by the explicit representation (\ref{2.12}), for every $(t,x)\in (0,T)\times {\mathbb R}^d$, we have
\begin{eqnarray*}
|u(t,x)|&\leq& \int_0^t\|f(s)\|_{L^p({\mathbb R}^d)}\|K(t-s)\|_{L^{\frac{p}{p-1}}({\mathbb R}^d)}dr\cr\cr&\leq& \|f\|_{L^\infty_q(0,t;L^p({\mathbb R}^d))}\int_0^ts^{-\frac{1}{q}}(t-s)^{-\frac{d}{2p}}ds\cr\cr&=&
t^{\frac{1}{2}}\|f\|_{L^\infty_q(0,t;L^p({\mathbb R}^d))}B(1-\frac{1}{q},\frac{1}{q}+\frac{1}{2}),
\end{eqnarray*}
where $B$ is the Beta function.

Therefore $u\in L^\infty(0,T;L^\infty({\mathbb R}^d))$ and
\begin{eqnarray}\label{2.13}
\|u\|_{L^\infty(0,T;L^\infty({\mathbb R}^d))}\leq CT^{\frac{1}{2}} \|f\|_{L^\infty_q(0,T;L^p({\mathbb R}^d))}.
\end{eqnarray}

For $x\in{\mathbb R}^d$ and $1\leq i\leq d$,
\begin{eqnarray}\label{2.14}
\Big|\partial_{x_i}u(t,x)\Big|&=&\Big|\int_0^t\int_{{\mathbb R}^d}\partial_{x_i}K(t-r,x-y)f(r,y)dydr\Big|\cr\cr
&\leq& \frac{1}{(2\pi)^{\frac{d}{2}}}\int_0^t\|f\|_{L^p({\mathbb R}^d)}(r)(t-r)^{-\frac{d}{2}-1}
\Big[\int_{{\mathbb R}^d}\Big|e^{-\frac{|x-y|^2}{2(t-r)}}|x_i-y_i|
\Big|^{\frac{p}{p-1}}dy\Big]^{\frac{p-1}{p}}dr\cr\cr
&\leq& C\|f\|_{L^\infty_q(0,t;L^p({\mathbb R}^d))}\int_0^tr^{-\frac{1}{q}}(t-r)^{-1+\frac{1}{q}}dr
\cr\cr&=&C(p,d) \|f\|_{L^\infty_q(0,T;L^p({\mathbb R}^d))},
\end{eqnarray}
where the constant in (\ref{2.14}) is given by
\begin{eqnarray}\label{2.15}
C(p,d)=\pi^{-\frac{d+p-1}{2p}}2^{\frac{p-d}{2p}}\Big[\Gamma\Big(\frac{2p-1}{2p-2}\Big)
\Big]^{\frac{p-1}{p}}\Big(\frac{p-1}{p}\Big)^{\frac{(d+1)p-d}{2p}}B(1-\frac{1}{q},\frac{1}{q}),
\end{eqnarray}
and $\Gamma$ is the gamma function. Since $1\leq i\leq d$ is arbitrary, $|\nabla u|\in L^\infty((0,T)\times{\mathbb R}^d)$.

Now we will show that $u\in \mathcal{C}([0,T];W^{1,p}({\mathbb R}^d))$. To prove this result, for $p\geq 1$, $\beta\geq 0$, let $H^{\beta,p}({\mathbb R}^d):=(I-\Delta)^{-\beta/2}(L^p({\mathbb R}^d))$ be the Bessel potential space with the norm
\begin{eqnarray*}
\|h\|_{H^{\beta,p}({\mathbb R}^d)}=\|(I-\Delta)^{-\beta/2}h\|_{L^p({\mathbb R}^d)}.
\end{eqnarray*}
For $h\in L^p({\mathbb R}^d)$, we use the notation $\mathcal{T}_th$ to denote $K(t,\cdot)\ast h$ with $K$ given in (\ref{2.12}). Then by a same discussion as \cite[Lemma 2.5]{Zha13}, we have the following claims:

(i) For $p>1$, $\beta>0$ and every $h\in L^p({\mathbb R}^d)$,  there is a constant $C(p,d,\beta)>0$ such that
\begin{eqnarray}\label{2.16}
\|\mathcal{T}_th\|_{H^{\beta,p}({\mathbb R}^d)}\leq C(p,d,\beta)t^{-\frac{\beta}{2}}\|h\|_{L^p({\mathbb R}^d)}.
\end{eqnarray}

(ii) For $p>1$, $\theta\in [0,1]$, there is a constant $C(p,d,\theta)>0$ such that for every $h\in H^{\beta,p}({\mathbb R}^d)$
\begin{eqnarray}\label{2.17}
\|\mathcal{T}_th-h\|_{L^p({\mathbb R}^d)}\leq C(p,d,\theta)t^{\frac{\theta}{2}}\|h\|_{H^{\theta,p}({\mathbb R}^d)}.
\end{eqnarray}

For every $0\leq s<t\leq T$, then
\begin{eqnarray}\label{2.18}
u(t,x)-u(s,x)&=&\int_0^t\mathcal{T}_{t-r}f(r)dr-\int_0^s\mathcal{T}_{s-r}f(r)dr\cr\cr&=&
\int_s^t\mathcal{T}_{t-r}f(r)dr+\int_0^s[\mathcal{T}_{t-r}-\mathcal{T}_{s-r}]f(r)dr\cr\cr&=&
\int_s^t\mathcal{T}_{t-r}f(r)dr+\int_0^s\mathcal{T}_{\frac{s-r}{2}}[\mathcal{T}_{t-s}-I]\mathcal{T}_{\frac{s-r}{2}}f(r)dr.
\end{eqnarray}
Observing that $W^{1,p}({\mathbb R}^d)=H^{1,p}({\mathbb R}^d)$ (see \cite{AF}), from (\ref{2.18}), (\ref{2.16}) and (\ref{2.17}), then
\begin{eqnarray}\label{2.19}
&&\|u(t)-u(s)\|_{W^{1,p}({\mathbb R}^d)}\cr\cr&\leq&
\int_s^t\|\mathcal{T}_{t-r}f(r)\|_{H^{1,p}({\mathbb R}^d)}dr+\int_0^s\|\mathcal{T}_{\frac{s-r}{2}}
[\mathcal{T}_{t-s}-I]\mathcal{T}_{\frac{s-r}{2}}f(r)\|_{H^{1,p}({\mathbb R}^d)}dr\cr\cr&\leq&
C\int_s^t(t-r)^{-\frac{1}{2}}\|f(r)\|_{L^p({\mathbb R}^d)}dr+C\int_0^s(s-r)^{-\frac{1}{2}}\|
[\mathcal{T}_{t-s}-I]\mathcal{T}_{\frac{s-r}{2}}f(r)\|_{L^p({\mathbb R}^d)}dr\cr\cr&\leq&
C\int_s^t(t-r)^{-\frac{1}{2}}\|f(r)\|_{L^p({\mathbb R}^d)}dr+C(t-s)^{\frac{\theta}{2}}\int_0^s(s-r)^{-\frac{1}{2}}
\|\mathcal{T}_{\frac{s-r}{2}}f(r)\|_{H^{\theta,p}({\mathbb R}^d)}dr\cr\cr&\leq&C\int_s^t
(t-r)^{-\frac{1}{2}}\|f(r)\|_{L^p({\mathbb R}^d)}dr+C(t-s)^{\frac{\theta}{2}}
\int_0^s(s-r)^{-\frac{1+\theta}{2}}
\|f(r)\|_{L^p({\mathbb R}^d)}dr,
\end{eqnarray}
for some $\theta\in [0,1]$.

Since $f\in L^\infty_q(0,T;L^p({\mathbb R}^d))$, if one chooses $\theta=(2(q-1)(q-2))/((3q-2)q)$, from (\ref{2.19}) by using H\"{o}lder's inequality, it then yields the following
\begin{eqnarray}\label{2.20}
&&\|u(t)-u(s)\|_{W^{1,p}({\mathbb R}^d)}\cr\cr&\leq&C\|f\|_{L^\infty_q(0,T;L^p({\mathbb R}^d))}
\Big[\int_s^t
(t-r)^{-\frac{1}{2}}r^{-\frac{1}{q}}dr+(t-s)^{\frac{\theta}{2}}
\int_0^s(s-r)^{-\frac{1+\theta}{2}}r^{-\frac{1}{q}}dr \Big ] \cr\cr&\leq& C\|f\|_{L^\infty_q(0,T;L^p({\mathbb R}^d))}|t-s|^{\frac{\theta}{2}}.
\end{eqnarray}
From this, one completes the proof for $g=0$.

For general $g$, since $u\in L^\infty(0,T;W^{1,\infty}({\mathbb R}^d))$, we conclude that: if $g\in L^\infty_q(0,T;L^p({\mathbb R}^d;{\mathbb R}^d))$, then $g\cdot\nabla u \in L^\infty_q(0,T;L^p({\mathbb R}^d))$. We define a
mapping from $L^\infty(0,T;W^{1,\infty}({\mathbb R}^d))$ to itself by
\begin{eqnarray*}
Tv(t,x)=\int_0^{t}K(t-s,\cdot)\ast (g(s,\cdot)\cdot \nabla v(s,\cdot))(x)ds+\int_0^{t}(K(t-s,\cdot)\ast f(s,\cdot))(x)ds.
\end{eqnarray*}
Observing that $\|g\|_{L^\infty_q(0,T;L^p({\mathbb R}^d))}$ is small enough, the mapping is contractive, so there is a unique $u\in L^\infty(0,T;W^{1,\infty}({\mathbb R}^d)$ satisfying $u=Tu$. This fact combining
an argument as $g=0$ implies the existence of generalized solutions of the Cauchy problem (\ref{2.8}).

Since (\ref{2.10}) holds, by virtue of (\ref{2.13}) and (\ref{2.14}), there is a constant $C_0(p,d)$ such that
\begin{eqnarray*}
\|u\|_{L^\infty(0,T;W^{1,\infty}({\mathbb R}^d))} \leq C_0(p,d)\Big[\|g\|_{L^\infty_q(0,T;L^p({\mathbb R}^d))}\|u\|_{L^\infty(0,T;W^{1,\infty}({\mathbb R}^d))}+\|f\|_{L^\infty_q(0,T;L^p({\mathbb R}^d))}\Big],
\end{eqnarray*}
which suggests that (\ref{2.11}) is valid since $\|g\|_{L^\infty_q(0,T;L^p({\mathbb R}^d))}$ is sufficiently small. $\Box$

\vskip4mm\noindent
\textbf{Remark 2.2.} (i) If $p,q\in (1,\infty)$, $f\in L^q(0,T;L^p({\mathbb R}^d))$ and $g\in L^q(0,T;L^p({\mathbb R}^d;{\mathbb R}^d))$, from classical $L^q(L^p)$ theory for second order parabolic PDEs, there is a unique $u\in W^{1,q}(0,T;L^p({\mathbb R}^d))\cap L^q(0,T;W^{2,p}({\mathbb R}^d))$ solving the Cauchy problem (\ref{2.8}). Using the Sobolev embedding theorems (or see \cite[Lemma 10.2]{KR}), if  $2/q+d/p<1$, then $\nabla u$ is bounded. But this embedding is not true in general when $2/q+d/p=1$. By extending the Banach space $L^q(0,T;L^p)$ to $L^\infty_q(0,T;L^p)$, under the critical case $2/q+d/p=1$, we also get the boundedness for $\nabla u$. In particular if we take $\tilde{f}\in L^\infty([0,T];L^p({\mathbb R}^d))$, $\tilde{g}\in L^\infty([0,T];L^p({\mathbb R}^d;{\mathbb R}^d))$, and set
\begin{eqnarray*}
f(t,x)=t^{-\frac{1}{q}}|\log t|^{-1}\tilde{f}(t,x), \  g(t,x)=t^{-\frac{1}{q}}|\log t|^{-1}\tilde{g}(t,x),
\end{eqnarray*}
then $f$ and $g$ are in $L^q(0,T;L^p)$. In this sense, we generalize the classical PDE's results.

(ii) When $p\geq 2$, one also proves $\partial_{x_i}u\in \mathcal{C}([0,T];L^p({\mathbb R}^d))$ for every $1\leq i\leq d$, by a duality method. Without loss of generality, here, we concentrate our attention for $g=0$. In this setting, we first recall that $L^\infty_q(0,T;L^p({\mathbb R}^d))\subset L^{q-}(0,T;L^p({\mathbb R}^d))$ ($=\cup_{1\leq r<q}L^r(0,T;L^p({\mathbb R}^d))$), with the help of \cite[Theorem 10.3]{KR}, $u$ given by (\ref{2.12}) is in $W^{1,q-}(0,T;L^p({\mathbb R}^d))\cap L^{q-}(0,T;W^{2,p}({\mathbb R}^d))$.

Let $\varrho_\varepsilon$ be a regularizing kernel on ${\mathbb R}$, i.e.
\begin{eqnarray}\label{2.21}
\varrho_\varepsilon =\frac{1}{\varepsilon} \varrho(\frac{\cdot}{\varepsilon}) \ \ with \ \ 0\leq \varrho \in \mathcal{C}^\infty_0({\mathbb R}) , \ \ support(\varrho)\subset (-1,1), \ \ \int_{{\mathbb R}}\varrho(r)dr=1.
\end{eqnarray}
We extend $u$ to the large interval $[-\varsigma,T+\varsigma]$ for $\varsigma>0$ and set $u_\varepsilon(t,x)=(u(\cdot,x)\ast \varrho_\varepsilon)(t)$. For $\varepsilon_1,\varepsilon_2>0$ and every $1\leq i\leq d$, we have
\begin{eqnarray}\label{2.22}
&&\frac{d}{dt}\|\partial_{x_i}u_{\varepsilon_1}(t)-\partial_{x_i}u_{\varepsilon_2}
(t)\|^p_{L^p({\mathbb R}^d)}\cr\cr&=& p\langle \partial_{x_i}u_{\varepsilon_1}^\prime(t)-\partial_{x_i}u_{\varepsilon_2}^\prime(t), [\partial_{x_i}u_{\varepsilon_1}(t)-\partial_{x_i}u_{\varepsilon_2}(t)]|\partial_{x_i}
u_{\varepsilon_1}(t)-\partial_{x_i}u_{\varepsilon_2}(t)|^{p-2}\rangle_{L^2({\mathbb R}^d)}
\cr\cr&=& -C(p)\langle u_{\varepsilon_1}^\prime(t)-u_{\varepsilon_2}^\prime(t), [\partial^2_{x_i}u_{\varepsilon_1}(t)-\partial^2_{x_i}u_{\varepsilon_2}(t)]|\partial_{x_i}
u_{\varepsilon_1}(t)-\partial_{x_i}u_{\varepsilon_2}(t)|^{p-2}\rangle_{L^2({\mathbb R}^d)},
\end{eqnarray}
where the prime in (\ref{2.22}) denotes the derivative in $t$ and $C(p)=p(p-1)$ in the last line.

From (\ref{2.22}), for every $0\leq s,t\leq T$, then
\begin{eqnarray*}
&&\|\partial_{x_i}u_{\varepsilon_1}(t)-\partial_{x_i}u_{\varepsilon_2}
(t)\|^p_{L^p({\mathbb R}^d)}\cr\cr&=&\|\partial_{x_i}u_{\varepsilon_1}(s)-\partial_{x_i}u_{\varepsilon_2}
(s)\|^p_{L^p({\mathbb R}^d)} \cr\cr&& -p(p-1)\int_s^t\langle u_{\varepsilon_1}^\prime(r)-u_{\varepsilon_2}^\prime(r), [\partial^2_{x_i}u_{\varepsilon_1}(r)-\partial^2_{x_i}u_{\varepsilon_2}(r)]|\partial_{x_i}
u_{\varepsilon_1}(r)-\partial_{x_i}u_{\varepsilon_2}(r)|^{p-2}\rangle_{L^2({\mathbb R}^d)}dr.
\end{eqnarray*}
Observing that for every $s\in (0,T)$, $\partial_{x_i}u_{\varepsilon_1}(s)\rightarrow \partial_{x_i}u(s)$ in $L^p({\mathbb R}^d)$ as $\varepsilon_1\rightarrow 0$, thus
\begin{eqnarray*}
&&\limsup_{\varepsilon_1,\varepsilon_2 \rightarrow 0}\sup_{0\leq t\leq T}\|\partial_{x_i}u_{\varepsilon_1}(t)-\partial_{x_i}u_{\varepsilon_2}
(t)\|^p_{L^p({\mathbb R}^d)}\cr\cr&\leq&C\limsup_{\varepsilon_1,\varepsilon_2 \rightarrow 0}\int_0^T|\langle u_{\varepsilon_1}^\prime(r)-u_{\varepsilon_2}^\prime(r), [\partial^2_{x_i}u_{\varepsilon_1}(r)-\partial^2_{x_i}u_{\varepsilon_2}(r)]|\partial_{x_i}
u_{\varepsilon_1}(r)-\partial_{x_i}u_{\varepsilon_2}(r)|^{p-2}\rangle_{L^2({\mathbb R}^d)}|dr
\cr\cr&\leq&C\limsup_{\varepsilon_1,\varepsilon_2 \rightarrow 0}\int_0^T\| u_{\varepsilon_1}^\prime-u_{\varepsilon_2}^\prime\|_{L^p({\mathbb R}^d)}\|\partial^2_{x_i}u_{\varepsilon_1}-\partial^2_{x_i}u_{\varepsilon_2}\|
_{L^p({\mathbb R}^d)}
\|\partial_{x_i}
u_{\varepsilon_1}-\partial_{x_i}u_{\varepsilon_2}\|_{L^p({\mathbb R}^d)}^{p-2}dr
\cr\cr&\leq&C\limsup_{\varepsilon_1,\varepsilon_2 \rightarrow 0}\int_0^T\| u_{\varepsilon_1}^\prime(r)-u_{\varepsilon_2}^\prime(r)\|_{L^p({\mathbb R}^d)}\|\partial^2_{x_i}u_{\varepsilon_1}(r)-\partial^2_{x_i}u_{\varepsilon_2}(r)\|
_{L^p({\mathbb R}^d)}dr
\cr\cr&\leq&C\limsup_{\varepsilon_1,\varepsilon_2 \rightarrow 0}\int_0^T\Big[\| u_{\varepsilon_1}^\prime(r)-
u_{\varepsilon_2}^\prime(r)\|_{L^p({\mathbb R}^d)}^2+
\|\partial^2_{x_i}u_{\varepsilon_1}(r)-\partial^2_{x_i}u_{\varepsilon_2}(r)\|^2
_{L^p({\mathbb R}^d)}\Big]dr=0.
\end{eqnarray*}
Thus $\{\partial_{x_i}u_{\varepsilon}\}_{0<\varepsilon<1}$ converge in $\mathcal{C}([0,T];L^p({\mathbb R}^d))$ to a function $\tilde{u}\in \mathcal{C}([0,T];L^p({\mathbb R}^d))$. Since we also know $\partial_{x_i}u_{\varepsilon}(t)$ converges to $\partial_{x_i}u(t)$ for a.e. $t\in [0,T]$, $u=\tilde{u}$.

(iii) From above proof, we claim that $u$ is continuous in $(t,x)$ and $\nabla u$ is bounded Borel measurable when $f\in L^\infty_q(0,T;L^p({\mathbb R}^d))$ and $g\in L^\infty_q(0,T;L^p({\mathbb R}^d;{\mathbb R}^d))$. If $f$ and $g$ are more regular, we also get the continuity of $\nabla u$.

\begin{corollary} \label{cor2.1} Suppose that $p$ and $q$ are given in Lemma \ref{lem2.1}, such that $f\in \mathcal{C}_q((0,T];L^p({\mathbb R}^d))$, $g\in \mathcal{C}_q((0,T];L^p({\mathbb R}^d;{\mathbb R}^d))$ and $\|g\|_{\mathcal{C}_q((0,T];L^p({\mathbb R}^d))}$ is small enough. Then there is a unique $u$ belonging to  $L^\infty(0,T;\mathcal{C}^1_b({\mathbb R}^d))\cap \mathcal{C}([0,T];W^{1,p}({\mathbb R}^d))$ solving the Cauchy problem (\ref{2.8}). Moreover,

(i) If $g=0$, $u\in L^\infty(0,T;\mathcal{C}^1_b({\mathbb R}^d))$.

(ii) If $f$ and $g$ meet in addition that $g\in \mathcal{C}^0_q((0,T];L^p({\mathbb R}^d;{\mathbb R}^d))$, $f\in \mathcal{C}^0_q((0,T];L^p({\mathbb R}^d))$, then $u\in \mathcal{C}([0,T];\mathcal{C}^1_b({\mathbb R}^d))$.
\end{corollary}
\vskip2mm\noindent
\textbf{Proof.} We only need to prove the continuity for $\nabla u$ and for simplicity, we show the case of $g=0$.

(i) For $x,h\in {\mathbb R}^d$ and $1\leq i\leq d$, by repeating the calculation in (\ref{2.14}), it yields that
\begin{eqnarray*}
\sup_{0\leq t\leq T}\Big|\partial_{x_i}u(t,x+h)-\partial_{x_i}u(t,x)\Big|
\leq C\sup_{0\leq t\leq T}\Big[\int_{{\mathbb R}^d}\Big|t^{\frac{1}{q}}f(t,x+h-z)-t^{\frac{1}{q}}f(t,x-z)\Big|^pdy\Big]^{\frac{1}{p}}.
\end{eqnarray*}
Notice that $t^{\frac{1}{q}}f\in \mathcal{C}([0,T];L^p({\mathbb R}^d))$, by an analogue argument of the proof of Proposition \ref{pro2.1}, if one lets $h$ approach to $0$, then
\begin{eqnarray*}
\sup_{0\leq t\leq T}\Big|\partial_{x_i}u(t,x+h)-\partial_{x_i}u(t,x)\Big|\rightarrow 0,
\end{eqnarray*}
which implies $\partial_{x_i}u\in L^\infty(0,T;\mathcal{C}({\mathbb R}^d))$. And $1\leq i\leq d$ is arbitrary, so $\nabla u$ is continuous in $x$ uniformly respect to $t$.

(ii) If the conditions $g\in \mathcal{C}^0_q((0,T];L^p({\mathbb R}^d;{\mathbb R}^d))$, $f\in \mathcal{C}^0_q((0,T];L^p({\mathbb R}^d))$  hold in addition, we will show the continuity of $\nabla u$ in $t$. We first show the continuity at $0$. In view of (\ref{2.14}), then for $t>0$,
\begin{eqnarray}\label{2.23}
\lim_{t\downarrow 0}\sup_{x\in {\mathbb R}^d}|\partial_{x_i}u(t,x)|\leq \lim_{t\downarrow 0}C[t^{\frac{1}{q}}\|f(t)\|_{L^p({\mathbb R}^d)}]=0.
\end{eqnarray}
Noting that $u(0,x)=0$, it implies $|\nabla u(0,x)|=0$, so $\nabla u$ is continuous in $t$ at $0$.

For $t>0$, we only prove the right continuous at $t$ since the proof for left continuous is the similar. Let $\vartheta>0$, $t\in (0,T)$ such that $t+\vartheta\in (0,T)$, then for every $1\leq i\leq d$
\begin{eqnarray*}
&&\Big|\partial_{x_i}u(t+\vartheta,x)-\partial_{x_i}u(t,x)\Big|\cr\cr&=&\Big|\int_0^{t+\vartheta}\int_{{\mathbb R}^d}\partial_{x_i}K(t+\vartheta-s,x-y)f(s,y)dyds-\int_0^t\int_{{\mathbb R}^d}\partial_{x_i}K(t-s,x-y)f(s,y)dyds\Big|\cr\cr
&\leq&
\Big|\int_0^\vartheta\int_{{\mathbb R}^d}\partial_{x_i}K(t+\vartheta-s,x-y)f(s,y)dyds\Big|\cr\cr&&+\Big|\int_0^t\int_{{\mathbb R}^d}\partial_{x_i}K(t-s,x-y)[f(s+\vartheta,y)-f(s,y)]dyds\Big|=:
I_1(t,\vartheta)+I_2(t,\vartheta).
\end{eqnarray*}
By using (\ref{2.23}), the $I_1(t,\vartheta)\rightarrow 0$ as $\vartheta\rightarrow 0$. Now let us calculate $I_2(t,\vartheta)$.
\begin{eqnarray*}
|I_2(t,\vartheta)|&\leq& C\int_0^t(t-s)^{-1+\frac{1}{q}}\|f(s+\vartheta)-f(s)\|_{L^p({\mathbb R}^d)}ds\cr\cr
&\leq&
C\int_0^t(t-s)^{-1+\frac{1}{q}}(s+\vartheta)^{-\frac{1}{q}}
\|(s+\vartheta)^{\frac{1}{q}}f(s+\vartheta)-s^{\frac{1}{q}}f(s)\|_{L^p({\mathbb R}^d)}ds
\cr\cr&&+C\int_0^t(t-s)^{-1+\frac{1}{q}}(s+\vartheta)^{-\frac{1}{q}}
\|(s+\vartheta)^{\frac{1}{q}}f(s)-s^{\frac{1}{q}}f(s)\|_{L^p({\mathbb R}^d)}ds
\cr\cr&\leq&
C\int_0^t(t-s)^{-1+\frac{1}{q}}s^{-\frac{1}{q}}
\|(s+\vartheta)^{\frac{1}{q}}f(s+\vartheta)-s^{\frac{1}{q}}f(s)\|_{L^p({\mathbb R}^d)}ds
\cr\cr&&+C\int_0^t(t-s)^{-1+\frac{1}{q}}s^{-\frac{1}{q}}
\|s^{\frac{1}{q}}f(s)\|_{L^p({\mathbb R}^d)}\Big[1-(\frac{s}{s+\vartheta})^{\frac{1}{q}}\Big]ds.
\end{eqnarray*}

Observing that as functions of $s$ on $(0,t)$, $(t-s)^{-1+\frac{1}{q}}s^{-\frac{1}{q}}\in L^1(0,t)$,
$\|(s+\vartheta)^{\frac{1}{q}}f(s+\vartheta)-s^{\frac{1}{q}}f(s)\|_{L^p({\mathbb R}^d)}$ and $\|s^{\frac{1}{q}}f(s)\|_{L^p({\mathbb R}^d)}
\Big[1-(s/(s+\vartheta))^{\frac{1}{q}}\Big]$ are bounded. By applying the dominated convergence theorem, then $I_2(t,\vartheta)\rightarrow 0$ as $\vartheta\rightarrow 0$. So we finish the proof. $\Box$

\vskip2mm\noindent
\textbf{Remark 2.3.} Even though a function $f$ lies in $\mathcal{C}^0_q((0,T];L^p({\mathbb R}^d))$, it is not in $L^q(0,T;L^p({\mathbb R}^d))$ generally. For example, if $\|f(t)\|_{L^p({\mathbb R}^d)}\leq C t^{-\frac{1}{q}}|\log t|^{-\beta}$ with some real number $\beta >0$ near $0$, then $f\in \mathcal{C}^0_q((0,T];L^p({\mathbb R}^d))$, and when $\beta>1/q$, it belongs to  $L^q(0,T;L^p({\mathbb R}^d))$, otherwise it does not lie in $L^q(0,T;L^p({\mathbb R}^d))$.

\vskip2mm\par
Let $W_t$ be a $d$-dimensional standard Wiener process, $X_0\in \mathcal{F}_0$, $\{\xi_t\}_{0\leq t\leq T}$ is a $\{\mathcal{F}_t\}_{0\leq t\leq T}$ adapted process, we define
\begin{eqnarray}\label{2.24}
X_t=X_0+\int_0^t\xi_sds+W_t.
\end{eqnarray}

We now give a Krylov type estimate.

\begin{theorem} \label{the2.1} Suppose $X_t$ is given by (\ref{2.24}) and $X_0=x\in {\mathbb R}^d$. Let $p,q\in [1,\infty)$ and $\mathcal{I}_Tf \in L^\infty_q(0,T;L^p({\mathbb R}^d))$ such that (\ref{1.4}) holds.  Let $C_0(p,d)$ be given in Lemma \ref{lem2.1}, then
\begin{eqnarray}\label{2.25}
{\mathbb E}\int_0^Tf(t,X_t)dt\leq C_0(p,d)\Big(1+{\mathbb E}\int_0^T|\xi_t|dt\Big)\|\mathcal{I}_Tf\|_{L^\infty_q(0,T;L^p({\mathbb R}^d))}.
\end{eqnarray}
\end{theorem}

\vskip2mm\noindent\textbf{Proof.} Let $u$ be given by
\begin{eqnarray*}
u(t,x)=\int_0^{t}(K(t-s,\cdot)\ast \mathcal{I}_Tf(s,\cdot))(x)ds.
\end{eqnarray*}
Observing that $\mathcal{I}_Tf\in L^\infty_q(0,T;L^p({\mathbb R}^d))$ with $p,q$ satisfying (\ref{1.4}), by virtue of Lemma \ref{lem2.1}, then $u\in \mathcal{C}([0,T];W^{1,p}({\mathbb R}^d))\cap L^\infty(0,T;W^{1,\infty}({\mathbb R}^d))$ and it solves the following Cauchy problem
\begin{eqnarray*}
\left\{\begin{array}{ll}
\partial_{t}u(t,x)=\frac{1}{2}\Delta u(t,x)+\mathcal{I}_Tf(t,x), \ (t,x)\in (0,T)\times {\mathbb R}^d, \\
u(0,x)=0, \  x\in{\mathbb R}^d,  \end{array}\right.
\end{eqnarray*}
in the sense of (\ref{2.9}). Moreover, (\ref{2.11}) holds with $g=0$.

By Sobolev embedding theorem, $\mathcal{C}([0,T];W^{1,p}({\mathbb R}^d))\subset \mathcal{C}([0,T]\times {\mathbb R}^d)$, $u$ is continuous in $t$ and $x$. If we smooth $u$ by convolution, then the modified function converges to $u$ for every $(t,x)\in [0,T]\times {\mathbb R}^d$, so we only concentrate our attention on $u\in \mathcal{C}^\infty([0,T]\times {\mathbb R}^d)$.

By virtue of It\^{o}'s formula, we have
\begin{eqnarray}\label{2.26}
&&du(T-t,X_t)\cr\cr&=&-\partial_t u(T-t,X_t)dt+\frac{1}{2}\Delta u(T-t,X_t)dt+\xi_t\cdot \nabla u(T-t,X_t)dt+\nabla u(T-t,X_t)\cdot dW_t
\cr\cr&=&\xi_t \cdot\nabla u(T-t,X_t)dt
-f(t,X_t)dt +\nabla u(T-t,X_t)\cdot dW_t.
\end{eqnarray}
Since $\nabla u$ is bounded, the last term in (\ref{2.26}) is a martingale, which implies that
\begin{eqnarray*}
u(T,x)={\mathbb E}\int_0^Tf(t,X_t)dt-{\mathbb E}\int_0^T\xi_t \cdot\nabla u(T-t,X_t)dt.
\end{eqnarray*}
Hence
\begin{eqnarray}\label{2.27}
{\mathbb E}\int_0^Tf(t,X_t)dt\leq \sup_{(t,x)\in (0,T)\times {\mathbb R}^d}|u(t,x)|+\sup_{(t,x)\in (0,T)\times {\mathbb R}^d}|\nabla u(t,x)|{\mathbb E}\int_0^T|\xi_t|dt.
\end{eqnarray}
By using Lemma \ref{lem2.1}, from (\ref{2.27}), (\ref{2.25})
holds true. $\Box$

\vskip2mm\noindent
\textbf{Remark 2.4.} Krylov's estimates will play a crucial role in proving the existence of weak solutions for SDE (\ref{1.1}). Observing that, when verifying a Krylov type estimate, the central part is to estimate the boundedness of $\nabla u$ ($u$ is the unique solution of a second order parabolic PDE). By finding this fact, we only assume that $f\in L^\infty_q(0,T;L^p({\mathbb R}^d))$ with $2/q+d/p=1$.

\section{Stochastic differential equations with irregular drift coefficients: the existence result }\label{sec3}
\setcounter{equation}{0}
We now consider the SDE (\ref{1.1}) and the main result is concerned with existence, which is stated below
\begin{theorem} \label{the3.1} Assume that $p,q\in [1,\infty)$. Let $b=b_1+b_2$ such that $\mathcal{I}_Tb_1 \in \mathcal{C}_q((0,T];L^p({\mathbb R}^d))$ with $p,q$ satisfying (\ref{1.4}), $b_2$ is bounded and Borel measurable. Suppose the constant $C_0(p,d)$ is given in (\ref{2.25}) and
\begin{eqnarray}\label{3.1}
\|\mathcal{I}_Tb_1\|_{\mathcal{C}_q((0,T];L^p({\mathbb R}^d))}<(2C_0(p,d))^{-1}.
\end{eqnarray}
There is a filtered probability space $(\tilde{\Omega},\tilde{\mathcal{F}},\{\tilde{\mathcal{F}}_t\}_{0\leq t\leq T},\tilde{{\mathbb P}})$, two  processes $\tilde{X}_t(x)$ and $\tilde{W}_t$ defined for $[0,T]$ on it such that $\tilde{W}_t$ is a $d$-dimensional $\{\tilde{\mathcal{F}}_t\}$-Wiener process and $\tilde{X}_t$ is an $\{\tilde{\mathcal{F}}_t\}$-adapted, continuous, $d$-dimensional process for which (\ref{1.5}) holds,
and almost surely, for all $t\in [0,T]$, (\ref{1.6}) holds.
\end{theorem}

\vskip2mm\noindent
\textbf{Proof.} We follow the proof of \cite[Theorem 1, p.87]{Kry80}. Firstly, we smooth out $b_i$ ($i=1,2$) using the convolution: $b_1^n(t,x)=(b_1(t,\cdot)\ast \rho_n)(x)$, $b_2^n(t,x)=(b_2(t,\cdot)\ast \rho_n)(x)$ with $\rho_n$ given by (\ref{2.2}).

According to Proposition \ref{pro2.1} and the properties of convolution, it is clear that, as $n\rightarrow \infty$,
\begin{eqnarray}\label{3.2}
\|\mathcal{I}_Tb_1^n-\mathcal{I}_Tb_1\|_{\mathcal{C}_q((0,T];L^p({\mathbb R}^d))} \rightarrow 0, \ \ \|b_2^n-b_2\|_{L^q(0,T;L^p_{loc}({\mathbb R}^d))} \rightarrow 0,
\end{eqnarray}
and for every $n\geq 1$,
\begin{eqnarray}\label{3.3}
\|\mathcal{I}_Tb_1^n\|_{\mathcal{C}_q((0,T];L^p({\mathbb R}^d))}\leq \|\mathcal{I}_Tb_1\|_{\mathcal{C}_q((0,T];L^p({\mathbb R}^d))}, \ \|b_2^n\|_{L^\infty((0,T)\times{\mathbb R}^d)}\leq \|b_2\|_{L^\infty((0,T)\times{\mathbb R}^d)}.
\end{eqnarray}
Moreover there is a sequence of integrable functions $h_i^n$ on $[0,T]$, such that
\begin{eqnarray*}
|b_i^n(t,x)-b_i^n(t,y)|\leq h_i^n(t)|x-y|, \ \ \forall \ x,y\in{\mathbb R}^d, \ \ i=1,2.
\end{eqnarray*}

Using Cauchy-Lipschitz's theorem, there is a unique $\{\mathcal{F}_t\}$-adapted, continuous, $d$-dimensional process $X_t^n(x)$ defined for $[0,T]$ on $(\Omega,\mathcal{F},\{\mathcal{F}_t\}_{0\leq t\leq T},{\mathbb P})$ such that
\begin{eqnarray}\label{3.4}
X_t^n=x+\int_0^tb^n(s,X^n_s)ds+W_t=x+\int_0^tb^n_1(s,X^n_s)ds+\int_0^tb^n_2(s,X^n_s)ds+W_t.
\end{eqnarray}

With the help of Theorem \ref{the2.1} and (\ref{3.3}),
\begin{eqnarray*}
{\mathbb E}\int_0^T|b_1^n(t,X^n_t)|dt&\leq&
\Big(1+
{\mathbb E}\int_0^T|b^n(t,X^n_t)|dt\Big)C_0(p,d)\|\mathcal{I}_Tb_1^n\|_{\mathcal{C}_q((0,T];L^p({\mathbb R}^d))}
\cr\cr&\leq&
\Big(1+T\|b_2\|_{L^\infty((0,T)\times{\mathbb R}^d)}
\cr\cr&&+{\mathbb E}\int_0^T|b_1^n(t,X^n_t)|dt\Big)
C_0(p,d)\|\mathcal{I}_Tb_1\|_{\mathcal{C}_q((0,T];L^p({\mathbb R}^d))}.
\end{eqnarray*}
Observing that (\ref{3.1}) holds, then
\begin{eqnarray*}
C_0(p,d)\|\mathcal{I}_Tb_1\|_{\mathcal{C}_q((0,T];L^p({\mathbb R}^d))}<\frac{1}{2},
\end{eqnarray*}
therefore
\begin{eqnarray}\label{3.5}
{\mathbb E}\int_0^T|b_1^n(t,X^n_t)|dt\leq
2\Big(1+T\|b_2\|_{L^\infty((0,T)\times{\mathbb R}^d)}\Big)C_0(p,d)\|\mathcal{I}_Tb_1\|_{\mathcal{C}_q((0,T];L^p({\mathbb R}^d))}.
\end{eqnarray}

On the other hand, $b_2^n$ is bounded uniformly in $n$, with the help of (\ref{3.3}), one concludes that
\begin{eqnarray}\label{3.6}
{\mathbb E}\int_0^T|b_2^n(t,X^n_t)|dt\leq T\|b_2\|_{L^\infty((0,T)\times {\mathbb R}^d)}.
\end{eqnarray}

By (\ref{3.4}) to (\ref{3.6}), then
\begin{eqnarray}\label{3.7}
\sup_{n}{\mathbb E}\int_0^T|X^n_t|dt\leq C<\infty.
\end{eqnarray}

If one replaces the time interval $(0,T)$ by $(t_1,t_2)$ for every $0\leq t_1<t_2\leq T$, similar calculations from (\ref{2.26}) to (\ref{2.27}) also yields that
\begin{eqnarray}\label{3.8}
&&{\mathbb E}\int_{t_1}^{t_2}|b_1^n(t,X^n_t)|dt\cr\cr&\leq& {\mathbb E}|U_n(T-t_2,X_{t_2}^n)-U_n(T-t_1,X_{t_1}^n)| +
 \sup_{(t,x)\in (0,T)\times {\mathbb R}^d}\|\nabla U_n\|{\mathbb E}\int_{t_1}^{t_2}|b^n(t,X^n_t)|dt
\cr\cr&\leq& \sup_{x\in{\mathbb R}^d}|U_n(T-t_2,x)-U_n(T-t_1,x)|\cr\cr&&+
\sup_{(t,x)\in (0,T)\times {\mathbb R}^d}\|\nabla U_n\|\Big({\mathbb E}|X_{t_2}^n-X_{t_1}^n|+{\mathbb E}\int_{t_1}^{t_2}|b^n(t,X^n_t)|dt\Big),
\end{eqnarray}
where $U_n$ is the unique generalized solution of
\begin{eqnarray*}
\left\{\begin{array}{ll}
\partial_{t}U_n(t,x)=\frac{1}{2}\Delta U_n(t,x)+|\mathcal{I}_Tb_1^n(t,x)|, \ (t,x)\in (0,T)\times {\mathbb R}^d, \\
U_n(0,x)=0, \  x\in{\mathbb R}^d,  \end{array}\right.
\end{eqnarray*}
With the aid of Sobolev's embedding theorem, (\ref{2.20}) and (\ref{3.3}), from (\ref{3.8})
\begin{eqnarray*}
&&{\mathbb E}\int_{t_1}^{t_2}|b_1^n(t,X^n_t)|dt\cr\cr&\leq&
C\|U_n(T-t_2)-U_n(T-t_1)\|_{W^{1,p}({\mathbb R}^d)}+
C_0(p,d)\|\mathcal{I}_Tb_1^n\|_{\mathcal{C}_q((0,T];L^p({\mathbb R}^d))}
\Big({\mathbb E}|X_{t_2}^n-X_{t_1}^n|\cr\cr&&+
|t_2-t_1|\|b_2^n\|_{L^\infty((0,T)\times{\mathbb R}^d)}
+{\mathbb E}\int_{t_1}^{t_2}|b_1^n(t,X^n_t)|dt\Big)
\cr\cr&\leq& C|t_2-t_1|^{\frac{\theta}{2}}\|\mathcal{I}_Tb_1^n\|_{\mathcal{C}_q((0,T];L^p({\mathbb R}^d))}+C_0(p,d)\|\mathcal{I}_Tb_1^n\|_{\mathcal{C}_q((0,T];L^p({\mathbb R}^d))}
\Big({\mathbb E}|X_{t_2}^n-X_{t_1}^n|\cr\cr&&+
|t_2-t_1|\|b_2^n\|_{L^\infty((0,T)\times{\mathbb R}^d)}
+{\mathbb E}\int_{t_1}^{t_2}|b_1^n(t,X^n_t)|dt\Big)
\cr\cr&\leq& C|t_2-t_1|^{\frac{\theta}{2}}\|\mathcal{I}_Tb_1\|_{\mathcal{C}_q((0,T];L^p({\mathbb R}^d))}+
C_0(p,d)\|\mathcal{I}_Tb_1\|_{\mathcal{C}_q((0,T];L^p({\mathbb R}^d))}
\Big({\mathbb E}|X_{t_2}^n-X_{t_1}^n|\cr\cr&&+
|t_2-t_1|\|b_2\|_{L^\infty((0,T)\times{\mathbb R}^d)}
+{\mathbb E}\int_{t_1}^{t_2}|b_1^n(t,X^n_t)|dt\Big),
\end{eqnarray*}
which suggests that
\begin{eqnarray}\label{3.9}
{\mathbb E}\int_{t_1}^{t_2}|b_1^n(t,X^n_t)|dt&\leq&
\frac{C_0(p,d)\|\mathcal{I}_Tb_1\|_{\mathcal{C}_q((0,T];L^p({\mathbb R}^d))}}
{1-C_0(p,d)\|\mathcal{I}_Tb_1\|_{\mathcal{C}_q((0,T];L^p({\mathbb R}^d))}}\Big[ C|t_2-t_1|^{\frac{\theta}{2}}\cr\cr&&+
{\mathbb E}|X_{t_2}^n-X_{t_1}^n|+
|t_2-t_1|\|b_2\|_{L^\infty((0,T)\times{\mathbb R}^d)}\Big],
\end{eqnarray}
where $\theta$ is given in (\ref{2.20}).

By (\ref{3.1}), there is a $\delta>0$, such that
$$
C_0(p,d)\|\mathcal{I}_Tb_1\|_{\mathcal{C}_q((0,T];L^p({\mathbb R}^d))}\leq (1-
C_0(p,d))\|\mathcal{I}_Tb_1\|_{\mathcal{C}_q((0,T];L^p({\mathbb R}^d))})(1-\delta).
$$
Combining (\ref{3.4}) and (\ref{3.9}), one reaches at
\begin{eqnarray*}
{\mathbb E}|X_{t_2}^n-X_{t_1}^n|&\leq& {\mathbb E}\int_{t_1}^{t_2}|b_1^n(t,X^n_t)|dt
+{\mathbb E}\int_{t_1}^{t_2}|b_2^n(t,X^n_t)|dt+{\mathbb E}|W_{t_2}-W_{t_1}|
\cr\cr&\leq& (1-\delta){\mathbb E}|X_{t_2}^n-X_{t_1}^n|+ C\Big(|t_2-t_1|^{\frac{\theta}{2}}+|t_2-t_1|+|t_2-t_1|^{\frac{1}{2}}\Big)
\cr\cr&\leq& (1-\delta){\mathbb E}|X_{t_2}^n-X_{t_1}^n|+ C|t_2-t_1|^{\frac{\theta}{2}},
\end{eqnarray*}
which implies that
\begin{eqnarray}\label{3.10}
{\mathbb E}|X_{t_2}^n-X_{t_1}^n|\leq \frac{C}{\delta}|t_2-t_1|^{\frac{\theta}{2}}\leq C|t_2-t_1|^{\frac{\theta}{2}}.
\end{eqnarray}

Combining (\ref{3.7}) and (\ref{3.10}) for every $\epsilon>0$, one concludes that
\begin{eqnarray}\label{3.11}
\lim_{c\rightarrow \infty}\sup_n\sup_{0\leq t\leq T}{\mathbb P}\{|X_t^n|>c\}=0,
\end{eqnarray}
and
\begin{eqnarray}\label{3.12}
\lim_{h\downarrow 0}\sup_n\sup_{0\leq t_1,t_2\leq T, |t_1-t_2|\leq h}{\mathbb P}\{|X_{t_1}^n-X_{t_2}^n|>\epsilon\}=0.
\end{eqnarray}

In view of Skorohod's representation theorem (see \cite[Lemma 2, p.87]{Kry80}), there is a probability space $(\tilde{\Omega},\tilde{\mathcal{F}},\tilde{{\mathbb P}})$ and random processes $(\tilde{X}_t^n,\tilde{W}_t^n)$, $(\tilde{X}_t,\tilde{W}_t)$ on this probability space such that

(i) the finite-dimensional distributions of $(\tilde{X}_t^n,\tilde{W}_t^n)$ coincide with the corresponding finite-dimensional distributions of the processes same as $(X_t^n,W_t^n)$.

(ii) $(\tilde{X}_\cdot^n,\tilde{W}_\cdot^n)$  converges to $(\tilde{X}_\cdot,\tilde{W}_\cdot)$, $\tilde{{\mathbb P}}$-almost surely.

In particular $\tilde{W}$ is still a Wiener process and
\begin{eqnarray}\label{3.13}
\tilde{X}_t^n=x+\int_0^tb_1^n(s,\tilde{X}^n_s)ds+\int_0^tb_2^n(s,\tilde{X}^n_s)ds+\tilde{W}^n_t.
\end{eqnarray}
For any $k\in\mN$, by virtue of Theorem \ref{the2.1},
\begin{eqnarray}\label{3.14}
&&\tilde{{\mathbb E}}\Big(\int_0^T|b_1^n(s,\tilde{X}^n_s)-b_1(s,\tilde{X}_s)|ds\Big)\cr\cr&\leq& \tilde{{\mathbb E}}\Big(\int_0^T|b_1^n(s,\tilde{X}^n_s)-b_1^k(s,\tilde{X}^n_s)|ds\Big)
+\tilde{{\mathbb E}}\Big(\int_0^T|b_1^k(s,\tilde{X}^n_s)-b_1^k(s,\tilde{X}_s)|ds\Big)
\cr\cr&&+\tilde{{\mathbb E}}\Big(\int_0^T|b_1^k(s,\tilde{X}_s)-b_1(s,\tilde{X}_s)|ds\Big)
\cr\cr&\leq& C\Big[ \|\mathcal{I}_Tb_1^n-\mathcal{I}_Tb_1^k\|_{\mathcal{C}_q((0,T];L^p({\mathbb R}^d))}+
\|\mathcal{I}_Tb_1^k-\mathcal{I}_Tb_1\|_{\mathcal{C}_q((0,T];L^p({\mathbb R}^d))}\Big]\cr\cr&&+
\tilde{{\mathbb E}}\Big(\int_0^T|b_1^k(s,\tilde{X}^n_s)-b_1^k(s,\tilde{X}_s)|ds\Big).
\end{eqnarray}
By letting $n\rightarrow \infty$ first, $k\rightarrow \infty$ next, from (\ref{3.2}) and (\ref{3.14}) we arrive at
\begin{eqnarray}\label{3.15}
\lim_{n\rightarrow \infty}\int_0^tb_1^n(s,\tilde{X}^n_s)ds=\int_0^tb_1(s,\tilde{X}_s)ds,\,\, \tilde{{\mathbb P}}-a.s..
\end{eqnarray}

Similarly, we obtain
\begin{eqnarray}\label{3.16}
&&\tilde{{\mathbb E}}\Big(\int_0^T|b_2^n(s,\tilde{X}^n_s)-b_2(s,\tilde{X}_s)|ds\Big)\cr\cr&\leq& \tilde{{\mathbb E}}\Big(\int_0^T|b_2^n(s,\tilde{X}^n_s)-b_2^k(s,\tilde{X}^n_s)|ds\Big)
+\tilde{{\mathbb E}}\Big(\int_0^T|b_2^k(s,\tilde{X}^n_s)-b_2^k(s,\tilde{X}_s)|ds\Big)
\cr\cr&&+\tilde{{\mathbb E}}\Big(\int_0^T|b_2^k(s,\tilde{X}_s)-b_2(s,\tilde{X}_s)|ds\Big)
\cr\cr&=&:J_1^n+J_2^n+J_3^n.
\end{eqnarray}
For $k$ fixed, as $n\rightarrow \infty$, $J_2^n\rightarrow 0$. For $J_1^n$, we have
\begin{eqnarray}\label{3.17}
J_1^n&\leq & \|b_2\|_{L^\infty((0,T)\times{\mathbb R}^d)}\tilde{{\mathbb E}}\int_0^T|1-1_{|\tilde{X}^n_s|\leq R}|ds +\tilde{{\mathbb E}}\Big(\int_0^T1_{|\tilde{X}^n_s|\leq R}|b_2^n(s,\tilde{X}^n_s)-b_2^k(s,\tilde{X}^n_s)|ds\Big)\cr\cr&\leq&
\frac{C}{R}\tilde{{\mathbb E}}\int_0^T|\tilde{X}^n_s|ds+C\|1_{|x|\leq R}|b_2^n-b_2^k|\|_{L^q(0,T;L^p({\mathbb R}^d))},
\end{eqnarray}
for every $R>0$.

By taking $n\rightarrow \infty$, $k\rightarrow \infty$, $R\rightarrow\infty$ in turn, then $J_1^n\rightarrow 0$ as $n\rightarrow \infty$. And the same conclusion for $J_3^n$ is true by a same discussion. Combining
(\ref{3.2}), (\ref{3.16}) and (\ref{3.17}) we arrive at
\begin{eqnarray}\label{3.18}
\lim_{n\rightarrow \infty}\int_0^tb_2^n(s,\tilde{X}^n_s)ds=\int_0^tb_2(s,\tilde{X}_s)ds,,\,\, \tilde{{\mathbb P}}-a.s..
\end{eqnarray}

From (\ref{3.13}), (\ref{3.15}) and (\ref{3.18}), one reaches at
\begin{eqnarray*}
\tilde{X}_t=x+\int_0^tb(s,\tilde{X}_s)ds+\tilde{W}_t.
\end{eqnarray*}
From this one finishes the proof since by using Theorem \ref{the2.1}, (\ref{3.2}) holds true obviously. $\Box$

Theorem \ref{the3.1} can be generalized to the case of nonconstant diffusion if the diffusion coefficients are regular enough. For simplicity, we give an application to the case of diffusion coefficient is time independent and $d=1$.
\begin{corollary} \label{cor3.1} Let $\sigma(x): {\mathbb R}\rightarrow {\mathbb R}$ be Borel measurable. Suppose that there are positive constants $\delta_1$ and $\delta_2$ such that $\delta_1\leq \sigma\leq \delta_2$. Consider the following SDE with nonconstant diffusion in ${\mathbb R}$
\begin{eqnarray}\label{3.19}
dX_t=b(t,X_t)dt+\sigma(X_t)dW_t, \ \ X_0=x\in {\mathbb R},  \ \ t\in [0,T].
\end{eqnarray}
Let $p$ and $q$ be given in Theorem \ref{the3.1}, that $b=b_1+b_2$ such that $\mathcal{I}_Tb_1 \in \mathcal{C}_q((0,T];L^p({\mathbb R}))$ and $\|\mathcal{I}_Tb_1\|_{\mathcal{C}_q((0,T];L^p({\mathbb R}))}$ is sufficiently small, $b_2$ is bounded Borel measurable.
Moreover, for this $p$, we assume in addition that $\sigma^\prime=\tilde{\sigma}_1+\tilde{\sigma}_2$, with $\tilde{\sigma}_1\in L^p({\mathbb R})$ and $\|\tilde{\sigma}_1\|_{L^p({\mathbb R})}$ is small enough, $\tilde{\sigma}_2\in L^\infty({\mathbb R})$. There is a filtered probability space $(\tilde{\Omega},\tilde{\mathcal{F}},\{\tilde{\mathcal{F}}_t\}_{0\leq t\leq T},\tilde{{\mathbb P}})$, two  processes $\tilde{X}_t(x)$ and $\tilde{W}_t$ defined for $[0,T]$ on it such that $\tilde{W}_t$ is a $1$-dimensional $\{\tilde{\mathcal{F}}_t\}$-Wiener process and $\tilde{X}_t$ is an $\{\tilde{\mathcal{F}}_t\}$-adapted, continuous, $1$-dimensional process for which (\ref{1.5}) holds, and almost surely, for all $t\in [0,T]$,
\begin{eqnarray*}
\tilde{X}_t=x+\int_0^tb(s,\tilde{X}_s)ds+\int_0^t\sigma(\tilde{X}_s)d\tilde{W}_s.
\end{eqnarray*}
\end{corollary}
\vskip2mm\noindent
\textbf{Proof.} The proof here is inspired by Zvonkin's transformation. Let us define
\begin{eqnarray}\label{3.20}
\Phi(x)=\int_0^x\frac{1}{\sigma(y)}dy,
\end{eqnarray}
and since $\delta_1\leq \sigma\leq \delta_2$, $\Phi^{-1}$ exists. Moreover, for every $x,y\in {\mathbb R}$,
\begin{eqnarray*}
\delta_2^{-1}|x-y|\leq |\Phi(x)-\Phi(y)|\leq \delta_1^{-1}|x-y|, \ \delta_1|x-y|\leq |\Phi^{-1}(x)-\Phi^{-1}(y)|\leq \delta_2|x-y|.
\end{eqnarray*}
Let us consider the following SDE
\begin{eqnarray}\label{3.21}
Y_t(y)=y+\int_0^t[b(s,\Phi^{-1}(Y_s))\sigma^{-1}(\Phi^{-1}(Y_s))-
\frac{1}{2}\sigma^\prime(\Phi^{-1}(Y_s))]ds+W_t.
\end{eqnarray}
Noting that $\mathcal{I}_Tb_1 \in \mathcal{C}_q((0,T];L^p({\mathbb R}))$ and $\|\mathcal{I}_Tb_1\|_{\mathcal{C}_q((0,T];L^p({\mathbb R}))}$ is small enough, so $\mathcal{I}_Tb_1(\cdot,\Phi^{-1}(\cdot))\in \mathcal{C}_q((0,T];L^p({\mathbb R}))$, $\|\mathcal{I}_Tb_1(\cdot,\Phi^{-1}(\cdot))\|_{\mathcal{C}_q((0,T];L^p({\mathbb R}))}$ is sufficiently small too. And this conclusion holds for $\tilde{\sigma}_1(\Phi^{-1})$. Upon using  Theorem \ref{the3.1}, there is a weak solution $(\tilde{Y}_t,\tilde{W}_t)$ of (\ref{3.21}) on a filtered probability space $(\tilde{\Omega},\tilde{\mathcal{F}},\{\tilde{\mathcal{F}}_t\}_{0\leq t\leq T},\tilde{{\mathbb P}})$ for which $\tilde{W}_t$ is a standard $1$-dimensional standard Wiener process to $\{\tilde{\mathcal{F}}_t\}_{0\leq t\leq T}$.

Initially, we assume $\sigma\in \mathcal{C}^\infty({\mathbb R})$, then $\Phi^{-1}$ is smooth. By utilising It\^{o}'s formula, one derives that
\begin{eqnarray}\label{3.22}
d\Phi^{-1}(\tilde{Y}_t)&=&[\Phi^{-1}]^\prime(\tilde{Y}_t)d\tilde{Y}_t-\frac{1}{2}
[\Phi^\prime(\Phi^{-1}(\tilde{Y}_t))]^{-1}\Phi^{\prime\prime}(\Phi^{-1}(\tilde{Y}_t))
[\Phi^{-1}]^\prime(\tilde{Y}_t)^2dt
\cr\cr&=&\sigma(\Phi^{-1}(\tilde{Y}_t))[b(t,\Phi^{-1}(\tilde{Y}_t))\sigma^{-1}(\Phi^{-1}(\tilde{Y}_t))
-\frac{1}{2}\sigma^\prime(\Phi^{-1}(\tilde{Y}_t))]dt
\cr\cr&&+\sigma(\Phi^{-1}(\tilde{Y}_t))d\tilde{W}_t+\frac{1}{2} \sigma(\Phi^{-1}(\tilde{Y}_t))\sigma^\prime(\Phi^{-1}(\tilde{Y}_t))dt\cr\cr&=&
b(t,\Phi^{-1}(\tilde{Y}_t))dt+\sigma(\Phi^{-1}(\tilde{Y}_t))d\tilde{W}_t,
\end{eqnarray}
which implies $(\tilde{X}_t,\tilde{W}_t)=(\Phi^{-1}(\tilde{Y}_t),\tilde{W}_t)$ is a weak solution of (\ref{3.19}).

For general $\sigma$, we smooth it by convolution $\sigma_\varepsilon=\sigma\ast\varrho_\varepsilon$ ($\varrho_\varepsilon$ is given in (\ref{2.21})). For $\Phi_\varepsilon^{-1}$, one gets an analogue of identity (\ref{3.22}). With the same argument as in Theorem \ref{the3.1}, by taking $\varepsilon\rightarrow 0$, one finishes the proof.  $\Box$

\vskip4mm\noindent
\textbf{Remark 3.1.} (i) The proof for the weak existence of solutions to SDE (\ref{3.19}) is inspired by Zvonkin's transformation. For more details in this topic, one consults to \cite{Zvo}.

(ii) When the Banach space $\mathcal{C}_q((0,T];L^p({\mathbb R}^d))$ is replaced by $L^q(0,T;L^p({\mathbb R}^d))$, with $q,p$ meeting condition (\ref{1.3}), there are many elegant study works \cite{Zha05,Zha11}. For example, under the hypothesises that

(1) $\sigma(t,x)$ is uniformly continuous in $x\in{\mathbb R}^d$ uniformly with respect to $t$ and there is a positive constant $\delta$ such that for all $(t,x)\in [0,T]\times{\mathbb R}^d$,
\begin{eqnarray*}
\delta|\xi|^2\leq\sum_{i,k}|\sigma^{i,k}(t,x)\xi_i|^2\leq \frac{1}{\delta}|\xi|^2.
\end{eqnarray*}

(2) $|\nabla_x\sigma_t(x)|, |b|\in L^q(0,T;L^p({\mathbb R}^d))$.
Zhang \cite{Zha11} obtained the existence and uniqueness of the strong solution to the following SDE
\begin{eqnarray}\label{3.23}
dX_t=b(t,X_t)dt+\sigma(t,X_t)dW_t, \ \ X_0=x\in {\mathbb R}^d,  \ \ t\in [0,T].
\end{eqnarray}

(iii) Other topics on SDE (\ref{3.23}) such as existence, uniqueness of solutions, stochastic homeomorphism, weak differentiability for $b$ and $\sigma$ in different classes, we refer the author to see \cite{Att,FF11,FGP2,GM01,MNP,YW} and the references cited up there.

\section{Uniqueness and the strong Feller property}\label{sec4}
\setcounter{equation}{0}
We first discuss the uniqueness. Before proceeding, let us present some useful lemmas.

\vskip2mm\par
Consider the SDE (\ref{3.23}), with $\sigma(t,x)\in{\mathbb R}^{d\times d}$. If $(X_t,W_t)$ is a weak solution on a probability space $(\Omega,\mathcal{F},{\mathbb P})$ with a reference family $\{\mathcal{F}_t\}_{0\leq t\leq T}$, for every $f\in \mathcal{C}^2_b({\mathbb R}^d)$, by It\^{o}'s formula we have
\begin{eqnarray*}
&&f(X_t)-f(x)-\int_0^tb(s,X_s)ds-\sum_{1\leq i,j\leq d}\int_0^ta_{i,j}(s,X_s)\partial_{x_i,x_j}^2f(X_s)ds\cr\cr&=&\sum_{1\leq i,k\leq d}\int_0^t\sigma_{i,k}(s,X_s)\partial_{x_i}f(X_s)dW_{k,s}.
\end{eqnarray*}
Hence
\begin{eqnarray}\label{4.1}
f(X_t)-f(x)-\int_0^tb(s,X_s)ds-\sum_{1\leq i,j\leq d}\int_0^ta_{i,j}(s,X_s)\partial_{x_i,x_j}^2f(X_s)ds\in \cM^c_2([0,T]),
\end{eqnarray}
if $\int_0^T|b(s,X_s)|ds<\infty$ and $\int_0^T|\sigma(s,X_s)|^2ds<\infty$, a.s., where $\cM^c_2([0,T])$ is the set of all continuous $\mathcal{F}_t$-adapted $L^2(0,T)$ martingale processes, $a_{i,j}=\sum_k\sigma_{i,k}\sigma_{j,k}/2$. Conversely, if a $d$-dimensional continuous adapted process $\{X_t\}_{0\leq t\leq T}$ defined on a probability space $(\Omega,\mathcal{F},{\mathbb P})$ with a reference family $\{\mathcal{F}_t\}_{0\leq t\leq T}$ satisfies (\ref{4.1}), then on an extension $(\tilde{\Omega},\tilde{\mathcal{F}},\tilde{{\mathbb P}})$ and $\{\tilde{\mathcal{F}}_t\}_{0\leq t\leq T}$, we can find a $d$-dimensional $\{\tilde{\mathcal{F}}_t\}_{0\leq t\leq T}$-Wiener process $\{\tilde{W}_t\}_{0\leq t\leq T}$ such that $(X,\tilde{W})$ is a weak solution of (\ref{3.23}) (see \cite[pp168-169]{IW}). And if $X$ meets (\ref{3.23}), its probability law ${\mathbb P}_x={\mathbb P}\circ X^{-1}$ on $d$-dimensional Wiener space $(W^d([0,T]),\mathcal{B}(W^d([0,T])))$ satisfies
\begin{eqnarray}\label{4.2}
f(w(t))-f(x)-\int_0^tb(s,w(s))ds-\sum_{1\leq i,j\leq d}\int_0^ta_{i,j}(s,w(s))\partial_{x_i,x_j}^2f(w(s))ds\in \cM^c_2,
\end{eqnarray}
for every $f\in \mathcal{C}_b^2({\mathbb R}^d)$. In summary, we have

\begin{lemma}(\cite[Proposition 2.1, p169]{IW}) \label{lem4.1} The existence of a weak solution of (\ref{3.23}) is equivalent to the existence of a $d$-dimensional process $X$ satisfying (\ref{4.1}), and this is also equivalent to the existence of a probability ${\mathbb P}$ on $(W^d([0,T]),\mathcal{B}(W^d([0,T])))$ satisfying (\ref{4.2}).
\end{lemma}

After then, we give another lemma, which is used to deal with the uniqueness
\begin{lemma}(\cite[Corollary, p206]{IW}) \label{lem4.2} If $(X,W)$ and $(X^\prime,W^\prime)$ are weak solutions (\ref{3.23}). Then ${\mathbb P}_x={\mathbb P}_x^\prime$ is equivalent to
\begin{eqnarray}\label{4.3}
\int_{W^d([0,T])}f(w(t)){\mathbb P}_x(dw)=\int_{W^d([0,T])}f(w(t)){\mathbb P}_x^\prime(dw),
\end{eqnarray}
for every $t\in [0,T]$ and every $f\in \mathcal{C}_b({\mathbb R}^d)$.
\end{lemma}

Now we give our uniqueness result.
\begin{theorem} \label{the4.1} Let $p,q$ and $b$ be stated in Theorem \ref{the3.1}. We assume that $\mathcal{I}_Tb_1\in \mathcal{C}^0_q((0,T];L^p({\mathbb R}^d;{\mathbb R}^d))$ in addition. If there are two weak solutions of (\ref{1.1}), then the probability laws of them on $d$-dimensional classical Wiener space $(W^d([0,T]),\mathcal{B}(W^d([0,T])))$ are the same.
\end{theorem}

\vskip2mm\noindent
\textbf{Proof.} We show the uniqueness by using It\^{o}-Tanack's trick (see \cite{FGP1}). Consider the following vector valued Cauchy problem
\begin{eqnarray}\label{4.4}
\left\{\begin{array}{ll}
\partial_tU(t,x)=\frac{1}{2}\Delta U(t,x)+\mathcal{I}_Tb_1(t,x)\cdot \nabla U(t,x)+\mathcal{I}_Tb_1(t,x), \ (t,x)\in (0,T)\times {\mathbb R}^d, \\
U(0,x)=0, \  x\in{\mathbb R}^d.  \end{array}\right.
\end{eqnarray}
According to Lemma \ref{lem2.1}, there is a unique generalized solution $U$. Moreover by Corollary \ref{cor2.1} (ii), $U\in \mathcal{C}([0,T];\mathcal{C}^1_b({\mathbb R}^d))$ and there is a $\delta>0$ such that
\begin{eqnarray*}
 \|U\|_{\mathcal{C}([0,T];\mathcal{C}^1_b({\mathbb R}^d))} \leq \frac{C_0(p,d)\|\mathcal{I}_T b_1\|_{\mathcal{C}^0_q((0,T];L^p({\mathbb R}^d))}}{1-C_0(p,d)
\|\mathcal{I}_Tb_1\|_{\mathcal{C}^0_q((0,T];L^p({\mathbb R}^d))}}<1-\delta,
\end{eqnarray*}
since (\ref{3.1}) holds.

We define $\Phi(t,x)=x+\mathcal{I}_TU(t,x)$, then it forms a diffeomorphism on ${\mathbb R}^d$, and
\begin{eqnarray}\label{4.5}
\delta<\|\nabla\Phi\|_{\mathcal{C}([0,T];\mathcal{C}_b({\mathbb R}^d))}< 2-\delta, \  \frac{1}{2-\delta}<\|\nabla\Psi\|_{\mathcal{C}([0,T];\mathcal{C}_b({\mathbb R}^d))}<\frac{1}{\delta},
\end{eqnarray}
where $\Psi(t,x)=\Phi^{-1}(t,x)$.

Using generalized It\^{o}'s formula (see \cite[Theorem 3.7]{KR}), if $(X_t,W_t)$ is a weak solution of (\ref{1.1}), then
\begin{eqnarray*}
 d\Phi(t,X_t(x))=b_2(t,X_t)dt+(I+\nabla \mathcal{I}_T U(t,X_t))dW_t, \ \Phi(0,X_0(x))=x+U(T,x).
\end{eqnarray*}
Denote $Y_t=X_t+\mathcal{I}_TU(t,X_t)$, then
\begin{eqnarray}\label{4.6}
 dY_t=b_2(t,\Psi(t,Y_t))dt+(I+\nabla \mathcal{I}_TU(t,\Psi(t,Y_t))dW_t, \ Y_0=y.
\end{eqnarray}

Now we assume that $(X,W)$ and $(X^\prime,W^\prime)$ are weak solutions of (\ref{1.1}) and the probability laws of $X$ and $X^\prime$ on $(W^d([0,T]),\mathcal{B}(W^d([0,T])))$
are ${\mathbb P}_x$ and ${\mathbb P}_x^\prime$ respectively. Correspondingly, we denote ${\mathbb P}_y$ and ${\mathbb P}_y^\prime$ the probability laws of $Y$ and $Y^\prime$. Since $Y_t=\Phi(t,X_t)$ and $\Phi\in \mathcal{C}([0,1];\mathcal{C}^1({\mathbb R}^d;{\mathbb R}^d))$ is a diffeomorphism on ${\mathbb R}^d$ uniformly for every $t\in[0,T]$, the relationships of ${\mathbb P}_x$ and ${\mathbb P}_y$, ${\mathbb P}_x^\prime$ and ${\mathbb P}_y^\prime$ are given by
\begin{eqnarray}\label{4.7}
{\mathbb P}_y={\mathbb P}_x\circ \Phi^{-1}, \ \  {\mathbb P}_y^\prime={\mathbb P}_x^\prime\circ \Phi^{-1},
\end{eqnarray}
where for a given measure $\mu$ on a Banach space $S$, and $\theta$ a map on $S$, we use the notation $\mu\circ \theta^{-1}$ to denote image measure of $\mu$ by the map $\theta$, i.e.
\begin{eqnarray*}
\int_S\phi d(\mu\circ \theta^{-1})=\int_S\phi(\theta) d\mu.
\end{eqnarray*}
Combining (\ref{4.5}) and noting that $b_2$ is bounded Borel measurable, $\mathcal{I}_T U$ is continuous in $(t,x)$, and $I+\mathcal{I}_TU$ meets uniformly elliptic condition, applying Lemma \ref{lem4.1} and \cite[Theorem 5.6]{SV1} (also see \cite[Theorem 3.3, p185]{IW} for time independent $\sigma$), the uniqueness of probability law for (\ref{4.6}) is true. So ${\mathbb P}_y={\mathbb P}_y^\prime$.

For every $f\in\mathcal{C}_b({\mathbb R}^d)$ and every $t\in [0,T]$, by  (\ref{4.7}) then \begin{eqnarray}\label{4.8}
\int_{W^d([0,T])}f(w(t)){\mathbb P}_x(dw)&=&\int_{W^d([0,T])}f(\Psi(t,\Phi(t,w(t)))){\mathbb P}_x(dw)
\cr\cr&=&\int_{W^d([0,T])}f(\Psi(t,w(t))){\mathbb P}_y(dw),
\end{eqnarray}
and
\begin{eqnarray}\label{4.9}
\int_{W^d([0,T])}f(w(t)){\mathbb P}_x^\prime(dw)&=&
\int_{W^d([0,T])}f(\Psi(t,\Phi(t,w(t)))){\mathbb P}_x^\prime(dw)
\cr\cr&=&\int_{W^d([0,T])}f(\Psi(t,w(t))){\mathbb P}_y^\prime(dw).
\end{eqnarray}
Since ${\mathbb P}_y={\mathbb P}_y^\prime$, and for every $t\in [0,T]$, $f\circ\Psi(t,\cdot)\in \mathcal{C}_b({\mathbb R}^d)$, from (\ref{4.8}) and (\ref{4.9}) one ends up with (\ref{4.3}). By applying Lemma \ref{lem4.2}, it is unique and we finish the proof. $\Box$

By using Zvonkin's transformation and in view of the fact the uniqueness in probability law implies the pathwise uniqueness in $d=1$ (see \cite[Proposition 1.1]{YTW}), we have

\begin{corollary} \label{cor3.2} Let $p,q$, $\sigma$ and $b$ satisfy the assumptions stated in Corollary \ref{cor3.1}. We suppose $\mathcal{I}_Tb_1 \in \mathcal{C}_q^0((0,T];L^p({\mathbb R}))$ further, then (\ref{3.19}) exists a unique strong solution.
\end{corollary}
\vskip2mm\noindent
\textbf{Proof.} Clearly, by Yamada, Watanabe'd theory (see \cite{YW}), it only needs to prove the uniqueness. Since the proof for uniqueness in probability law  is analogue of the proof of Theorem \ref{the4.1}, we prove the pathwise uniqueness only. On the other hand, by the relationship between (\ref{3.19}) and (\ref{3.21}), it suffices to show the pathwise for (\ref{1.1}) on $d=1$. When $b$ is time independent and, bounded and continuous in $x$, the proof can be seen in \cite[Lemma p.75]{YTW}. But the proof there is adapted for the present $b$, so we completes the proof. $\Box$

\vskip2mm\noindent
\textbf{Remark 4.1.} Here we do not argue the general case, i.e. $\sigma$ is time dependent and $d>1$. As discussed in \cite{Kry69}, we may consider the existence and uniqueness for weak solutions, such that the uniqueness holds only in the sense of finite dimensional probability laws. For more details in this topic one can refers to \cite{Kry69} and the references cited up there.

Now we discuss the Feller property and the existence of density and initially we give a lemma.
\begin{lemma}(\cite{SV2}) \label{lem4.3} Consider the following SDE
\begin{eqnarray}\label{4.10}
\left\{\begin{array}{ll}
dX_t(x)=b(t,X_t(x))dt+\sigma(t,X_t(x))dW_t, \ s<t\leq T,\\
X_s(x)=x\in{\mathbb R}^d,  \end{array}\right.
\end{eqnarray}
Suppose that $b$ is bounded and Borel measurable, $\sigma$ is bounded continuous and $(a_{i,j})=(\sum_k\sigma_{i,k}\sigma_{i,j})/2$ is uniformly elliptic. Then there is a unique weak solution of (\ref{4.10}) which is a strong Markov process. Let $\tilde{P}(s,x,t,dy)$ be the transition probabilities and for every bounded function $f$, we define
\begin{eqnarray}\label{4.11}
\tilde{P}_{s,t}f(x)=\int_{{\mathbb R}^d}f(y)\tilde{P}(s,x,t,dy),
\end{eqnarray}
then we have the following claims:

(i) $\tilde{P}_{s,t}f(x)$ is continuous in $s$ and $x$ for $s<t$.

(ii) $\tilde{P}(s,x,t,dy)$ has a density $\tilde{p}(s,x,t,y)$ for almost all $t\in [0,T]$ which satisfies
\begin{eqnarray}\label{4.12}
\int_{t_0}^T\int_{{\mathbb R}^d}|\tilde{p}(s,x,t,y)|^rdydt<\infty.
\end{eqnarray}
for every $r\in [1,\infty)$ provided $s<t_0$.
\end{lemma}

Now we give our second result.
\begin{theorem} \label{the4.2}  Let $p,q$ and $b$ be described in Theorem \ref{the4.1}. Consider SDE (\ref{4.10}) with $\sigma=I_{d\times d}$. Then there is a unique weak solution. Let $P(s,x,t,dy)$ be the transition probabilities and for every bounded function $f$, we set $P_{s,t}f(x)$ by (\ref{4.12}). Then $P_{s,t}f(x)$ is continuous in $s$ and $x$ for $s<t$ and $P(s,x,t,dy)$ has a density $p(s,x,t,y)$ for almost all $t\in [0,T]$ which meets (\ref{4.12}).
\end{theorem}
\noindent
\textbf{Proof.} By Theorems  \ref{the3.1} and \ref{the4.1}, there is a unique weak solution of (\ref{4.10}). Moreover, by (\ref{4.7}) and  Lemma \ref{lem4.3}, ${\mathbb P}_{s,x}$ is a strong Markov process. Let $P(s,x,t,dy)$ be its transition probabilities, then for every $f\in L^\infty({\mathbb R}^d)$,
\begin{eqnarray}\label{4.13}
P_{s,t}f(x)={\mathbb E}^{{\mathbb P}_{s,x}}f(w(t))=\int_{{\mathbb R}^d}f(y)P(s,x,t,dy).
\end{eqnarray}
Correspondingly, if one argues SDE (\ref{4.6}) with initial data from time $t=s$ and denotes the transition probabilities by $\tilde{P}(s,\Phi(s,x),t,dy)$, then with the help of (\ref{4.8}), it yields that
\begin{eqnarray}\label{4.14}
P_{s,t}f(x)=\tilde{P}_{s,t}f(\Psi(t,x))=\int_{{\mathbb R}^d}f(\Psi(t,y))\tilde{P}(s,\Phi(s,x),t,dy).
\end{eqnarray}
So $P_{s,t}f(x)$ is continuous in $s$ and $x$ for $s<t$. In particular, the semi-group $P_{0,t} (=:P_t)$ has strong Feller property. And by Lemma \ref{lem4.3}, $\tilde{P}_{s,t}$ has a density $\tilde{p}(s,\Phi(s,x),t,y)$, from (\ref{4.14}), then
\begin{eqnarray*}
P_{s,t}f(x)=\int_{{\mathbb R}^d}f(\Psi(t,y))
\tilde{p}(s,\Phi(s,x),t,y)dy=\int_{{\mathbb R}^d}f(y)
\tilde{p}(s,\Phi(s,x),t,\Phi(t,y))|\nabla \Phi(t,y)|dy.
\end{eqnarray*}
Hence, for almost all $t\in [0,T]$, $P_{s,t}$ has a density $p(s,x,t,y)$ which is given by
\begin{eqnarray}\label{4.15}
p(s,x,t,y)=\tilde{p}(s,\Phi(s,x),t,\Phi(t,y))|\nabla \Phi(t,y)|.
\end{eqnarray}
From (\ref{4.15}),for every $r\in [1,\infty)$ and $s<t_0$
\begin{eqnarray*}
\int_{t_0}^T\int_{{\mathbb R}^d}|p(s,x,t,y)|^rdydt&=&
\int_{t_0}^T\int_{{\mathbb R}^d}|\tilde{p}(s,\Phi(s,x),t,\Phi(t,y))|\nabla \Phi(t,y)||^rdydt\cr\cr&=&\int_{t_0}^T\int_{{\mathbb R}^d}|\tilde{p}(s,\Phi(s,x),t,y)|\nabla \Phi(t,\Psi(t,y))||^r|\nabla \Psi(t,y)|dydt
\cr\cr&\leq&C\int_{t_0}^T\int_{{\mathbb R}^d}|\tilde{p}(s,\Phi(s,x),t,y)|^rdydt<\infty,
\end{eqnarray*}
where in the last line we have used (\ref{4.5}). $\Box$

\section*{Acknowledgements}
This research was partly supported by the NSF of China grants 11501577, 11771123.

\end{document}